\pdfoutput=1 

\documentclass[twocolumn]{svjour3mod}

\usepackage{algorithm2e}
\usepackage{amsfonts}
\usepackage{amsmath}
\usepackage{amsopn}
\usepackage{amssymb}
\usepackage{blkarray}
\usepackage{graphicx}
\usepackage{lipsum}
\usepackage{multirow}
\usepackage{tikz}
\usepackage{tikz-cd}
\usepackage[customcolors]{hf-tikz}
\usepackage{xspace}

\usepackage[final]{listings}

\usepackage{epstopdf}

\usepackage[colorlinks]{hyperref}
\hypersetup{
  linkcolor = blue,
  citecolor = red,
  urlcolor  = brown
}

\hfsetfillcolor{blue!10}
\hfsetbordercolor{blue}
 
\DeclareMathOperator{\arctantwo}{arctan2}
\newcommand{\BLAS}{\textsc{BLAS}\xspace}
\newcommand{\code}[1]{\texttt{#1}}
\newcommand{\degrees}{\ensuremath{^\circ}}
\DeclareMathOperator*{\diagoper}{\text{diag}}
\newcommand{\ebox}[1]{\hbox{\ensuremath{#1}}}
\newcommand{\EigValMatrix}{\Lambda}
\newcommand{\EigVecMatrix}{V}
\newcommand{\Expectation}{\mathrm{E}}
\newcommand{\eigenvector}{v}
\newcommand{\Identity}{I}
\newcommand{\inverse}{^{-1}}
\newcommand{\kronvector}{\boldsymbol{\delta}}
\newcommand{\LAPACK}{\textsc{LAPACK}\xspace}
\newcommand{\ltwonorm}[1]{\Vert #1\Vert}
\newcommand{\modbesselfcxnfirst}[1]{I_{#1}}
\newcommand{\orientedeigvec}{\mathcal{V}}
\newcommand{\PCA}{\textsc{PCA}\xspace}
\newcommand{\realsspace}{\mathbb{R}}
\DeclareMathOperator{\sign}{sign}
\newcommand{\spherespace}{\mathbb{S}}
\newcommand{\SVD}{\textsc{SVD}\xspace}
\newcommand{\TensorFlow}{\textsc{TensorFlow}\xspace}
\newcommand{\transpose}{^T}
\newcommand{\vMF}{\textsc{vMF}\xspace}
\newcommand{\vmfconc}{\kappa}
\newcommand{\vmfmean}{\mu}
\newcommand{\vmfpdf}{p_{\text{vmf}}}
\newcommand{\VonMisesFisher}{\textsc{VonMisesFisher}\xspace}

\journalname{}

\begin{document}

\title{A Consistently Oriented Basis for Eigenanalysis}

\author{Jay Damask}
\authorrunning{Jay Damask}

\institute{
\at New York, NY. \\ 
\email{jaydamask@buell-lane-press.co} \\ 
orcid: 0000-0003-3288-2100 
}

\date{}

\maketitle

\begin{abstract}
%
%
%
%
Repeated application of machine-learning, eigen-centric methods to an evolving dataset reveals that eigenvectors calculated by well-established computer implementations are not stable along an evolving sequence. This is because the sign of any one eigenvector may point along either the positive or negative direction of its associated eigenaxis, and for any one eigen call the sign does not matter when calculating a solution. This work reports an algorithm that creates a consistently oriented basis of eigenvectors. The algorithm postprocesses any well-established eigen call and is therefore agnostic to the particular implementation of the latter. Once consistently oriented, directional statistics can be applied to the eigenvectors in order to track their motion and summarize their dispersion. When a consistently oriented eigensystem is applied to methods of machine-learning, the time series of training weights becomes interpretable in the context of the machine-learning model. Ordinary linear regression is used to demonstrate such interpretability. A reference implementation of the algorithm reported herein has been written in Python and is freely available, both as source code and through the \code{thucyd} Python package.

\end{abstract}


%
%

\section{Overview}
\label{sec: overview}

Eigenanalysis plays a central role in many aspects of machine learning because the natural coordinates and relative lengths of a multifactor dataset are revealed by the analysis \cite{Hastie:2001,Alpaydin:2010,Kuhn:2013,Scholkopf:2002}. Singular-value decomposition and eigendecomposition are two of the most common instances of eigenanalysis \cite{Strang:1988,Golub:2014}, yet there is a wide range of adaptation, such as eigenface \cite{FaceRecognition:2011}. The eigenanalysis ecosystem is so important and well studied that highly optimized implementations are available in core numerical libraries such as \LAPACK \cite{lapack:1999}. 

This article does not seek to suggest or make any improvement to the existing theoretical or numerical study of eigenanalysis as it applies to any one standalone problem. Instead, the focus of this article is on the improvement of eigenanalysis when applied repeatedly to an evolving dataset. With today's world being awash in data, and with machine-learning techniques being applied regularly, if not continuously, on ever-expanding datasets as new data arrives, the \emph{evolution} of the eigensystem ought to be included in any well-conceived study. Yet, the researcher will find that the \emph{orientation} of the eigenvector basis produced by current eigenanalysis algorithms is not consistent, that it may change according to hardware and software, and that it is sensitive to data perturbations. 

\begin{figure*}[t]\sidecaption
  \centering
  \includegraphics[width=126mm]{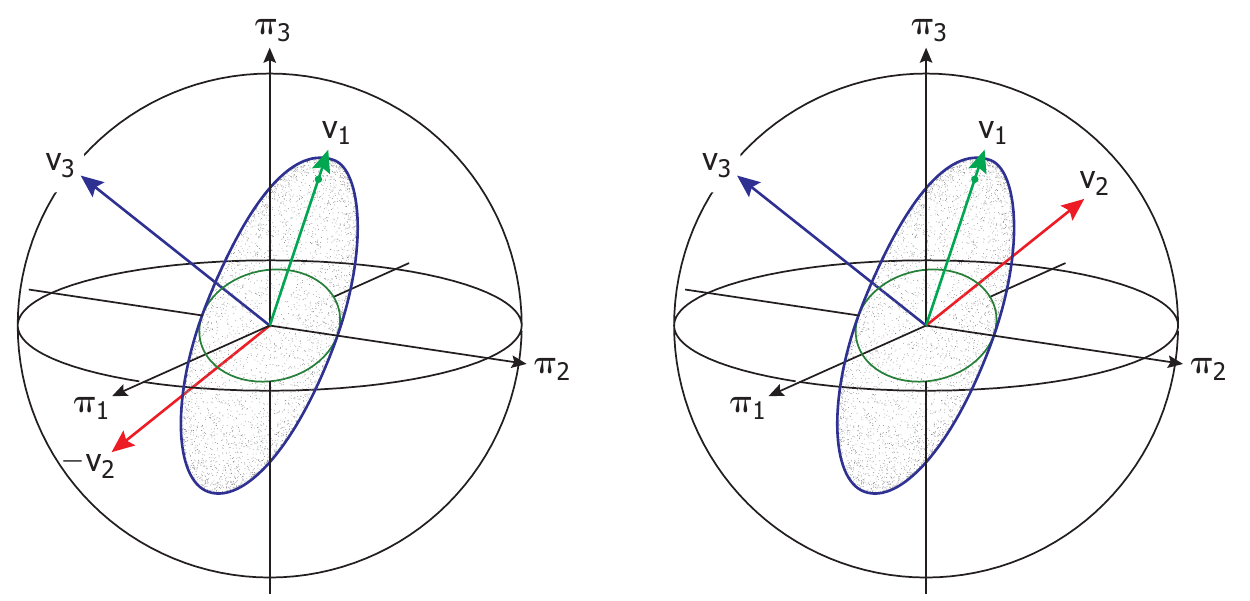}
  \caption{An ellipsoid-shaped point cloud in~$\realsspace^3$ typical of centered, variance-bound data having a stationary or slowly changing underlying random process. The point cloud is quoted in a constituent basis \ebox{(\pi_1, \pi_2, \pi_3)} while the eigenvectors lie along the natural axes of the data. Here, the ambiguity of the direction of~$v_2$ is unresolved. \emph{Left}, a left-handed basis of eigenvectors. \emph{Right}, a right-handed basis of eigenvectors.
  }
  \label{fig: r3 point clouds}
\end{figure*}

This paper demonstrates a method for the \emph{consistent orientation} of an eigensystem, in which the orientation is made consistent after an initial call to a canonical eigenanalysis routine. This method, however, does not address the issue of a \emph{correct} orientation: correctness is context specific and must be addressed in a manner exogenous to eigenanalysis. The analyst will need to know, for instance, whether a positive change in the factors aligned to particular eigenvectors should imply a positive or negative change in a dependent feature. In this article, consistency alone is sought.

As a consequence of the method detailed below, the evolution of the eigenvectors themselves can be tracked. Eigenvectors that span~$\realsspace^n$ lie on a hypersphere $\spherespace^{n-1} \subseteq \realsspace^n$, and the mean pointing direction on the sphere and dispersion about the mean are statistics of interest, as is any statistically significant drift. Basic results from directional statistics (see \cite{Mardia:2000,Jammal:2001,Ley:2017}) are used to study the mean direction and dispersion of the common and differential modes over eigenvector evolution.

%
%

\section{Inconsistent Orientation}
\label{sec: inconsistent orientation}

Given a square, symmetric matrix, $M\transpose =M$, with real-valued entries, $M \in \realsspace^{n\times n}$, the eigendecomposition is
\begin{equation}
  M = \EigVecMatrix \EigValMatrix \EigVecMatrix\transpose,
  \label{eq: eigen decomp def}
\end{equation}
where the columns of~$\EigVecMatrix \in \realsspace^{n\times n}$ are the eigenvectors of~$M$, the diagonal entries along~$\EigValMatrix \in \realsspace^{n\times n}$ are the eigenvalues, and entries in both are purely real. Moreover, when properly normalized, \ebox{\EigVecMatrix\transpose\EigVecMatrix = \Identity}, indicating that the eigenvectors form an orthonormal basis. 

To inspect how the eigenvectors might be oriented, write out~$\EigVecMatrix$ as \ebox{\left( \eigenvector_1, \eigenvector_2, \ldots, \eigenvector_n \right)}, where~\ebox{v_k} is a column vector, and take the inner product two ways:
\begin{equation*}
  \begin{aligned}
      \left( \eigenvector_1, \eigenvector_2, \ldots, \eigenvector_n \right)
      \transpose
      \left( \eigenvector_1, \eigenvector_2, \ldots, \eigenvector_n \right)  & = \Identity \quad\text{and} \\
      \left( \eigenvector_1, -\eigenvector_2, \ldots, \eigenvector_n \right)
      \transpose
      \left( \eigenvector_1, -\eigenvector_2, \ldots, \eigenvector_n \right) & = \Identity.
  \end{aligned}
\end{equation*}
A choice of~$+\eigenvector_2$ is indistinguishable from a choice of~$-\eigenvector_2$ because in either case the equality in~(\ref{eq: eigen decomp def}) holds. In practice, the eigenvectors take an \emph{inconsistent orientation}, even though they will have a particular orientation once returned from a call to an eigenanalysis routine.

A similar situation holds for singular-value decomposition (\SVD). For a centered panel of data~$P \in \realsspace^{m\times n}$ where \ebox{m>n}, \SVD gives
\begin{equation}
  P = U\: \EigValMatrix^{1/2}\: \EigVecMatrix\transpose,
  \label{eq: SVD decomp def}
\end{equation}
where the columns in panel~$U \in \realsspace^{m\times n}$ are the projections of the columns of~$P$ onto the eigenvectors~$\EigVecMatrix$, scaled to have unit variance. The eigenvectors have captured the natural orientation of the data. Unlike the example above, the eigensystem here is positive semidefinite, so all eigenvalues in~$\EigValMatrix$ are nonnegative, and in turn a scatter plot of data in~$P$ typically appears as an ellipsoid in~$\realsspace^n$ (see figure~\ref{fig: r3 point clouds}). Following the thought experiment above, a sign flip of an eigenvector in~$\EigVecMatrix$ flips the sign of the corresponding column in~$U$, thereby preserving the equality in~(\ref{eq: SVD decomp def}). Consequently, the eigenvector orientations produced by \SVD also have inconsistent orientation, since either sign satisfies the equations. 

In isolation, such inconsistent eigenvector orientation does not matter. The decomposition is correct and is generally produced efficiently. However, in subsequent analysis of the eigensystem---principal component analysis (\PCA) being a well known example \cite{Hastie:2001}---and for an evolving environment, the inconsistency leads to an inability to \emph{interpret} subsequent results and the inability to track the eigenvector evolution.

%
\begin{figure*}[t]\sidecaption
  \centering
  \includegraphics[width=126mm]{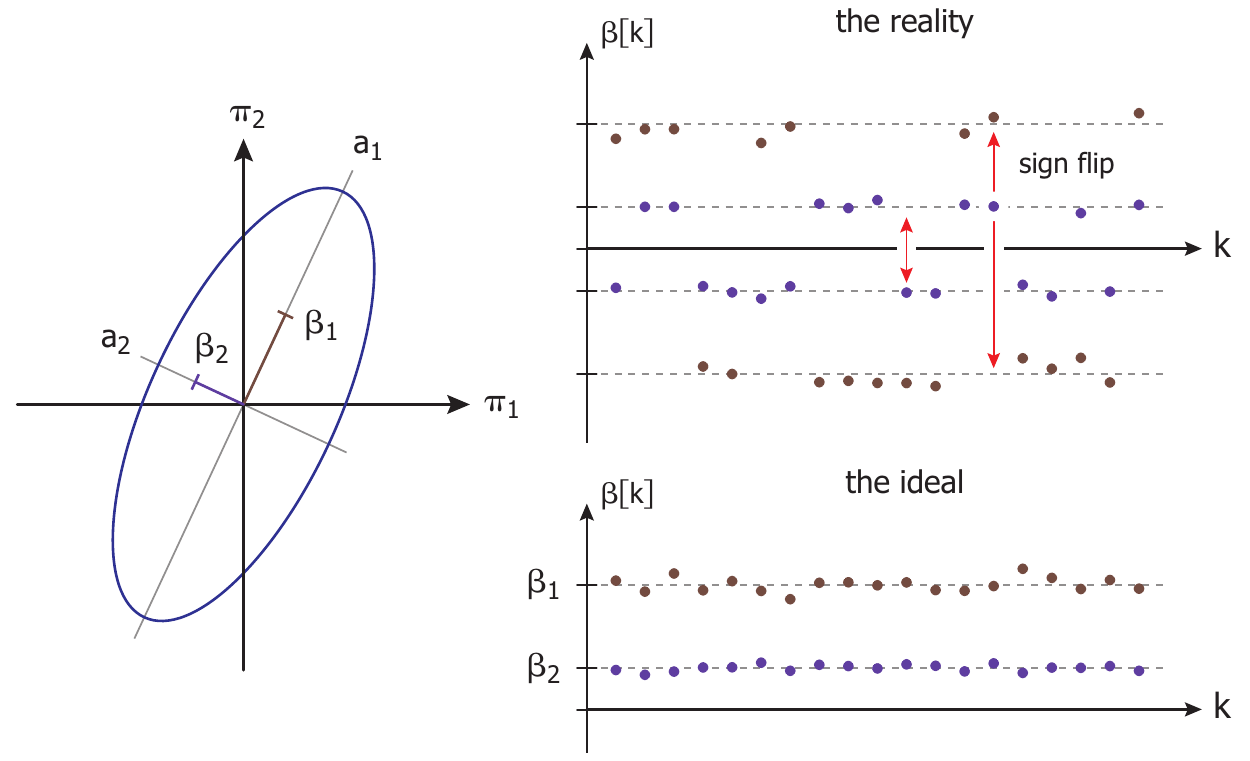}
  \caption{\emph{Left}, in $\realsspace^2$ a regression along the eigenaxes yields the tuple \ebox{(\beta_1, \beta_2)}. \emph{Right, top}, because the eigenvector direction can flip from one eigenanalysis to the next, the \emph{signs} of the $\beta$-tuple also flip, defeating any attempt at interpretability. \emph{Right, bottom}, ideally there is no sign flipping of the eigenvectors and therefore no flips of the $\beta$-tuple along a sequence. 
  }
  \label{fig: r2 beta sign flips}
\end{figure*}

For instance, a simple linear regression performed on the eigensystem (whether with reduced dimension or not) with dependent variable~$y \in \realsspace^m$ reads
\begin{equation}
  y = U \EigValMatrix^{1/2}\: \beta + \epsilon,
  \label{eq: OLS on eigensystem example}
\end{equation}
where the weight vector is~$\beta \in \realsspace^{n}$, \ebox{\beta = \left(\beta_1, \ldots, \beta_n\right)\transpose}. Each entry in~$\beta$ corresponds to a column in~$U$, which in turn corresponds to an eigenvector in~$\EigVecMatrix$. If the sign of one of the eigenvectors is flipped (compare \emph{left} and \emph{right} in figure~\ref{fig: r3 point clouds}), the corresponding column in~$U$ has its sign flipped, \emph{which flips the sign of the corresponding~$\beta$ entry}. In an evolving environment, a particular regression might be denoted by index~$k$, and the $k$th regression reads
\begin{equation}
  y_k = U_k \EigValMatrix_k^{1/2} \: \beta_k + \epsilon_k.
  \label{eq: OLS on eigensystem example 2}
\end{equation}
When tracking any particular entry in the~$\beta$ vector, say~$\beta_i$, it will generally be found that~$\beta_{i,k}$ flips sign throughout the evolution of system (see figure~\ref{fig: r2 beta sign flips}). There are two consequences: First, for a single regression, the \emph{interpretation} of how a change of the~$i$th factor corresponds to a change in~$y$ is unknown because the~$i$th eigenvector can take either sign. Second, as the data evolves, the sign of~$\beta_{i,k}$ will generally change, even if the unsigned $\beta_i$ value remains fixed, so that the time series of~$\beta$ cannot be meaningfully analyzed. If these sources of sign flips can be eliminated, then a time series of weights~$\beta$ can be interpreted as we might have expected (see figure~\ref{fig: r2 beta sign flips}~\emph{lower right}).

%
%

\section{Toward a Consistent Orientation}
\label{sec: toward a consistent orientation}

To construct a consistent eigenvector basis we must begin with a reference, and \SVD offers a suitable start. The data in matrix~$P$ of~(\ref{eq: SVD decomp def}) is presented with~$n$ columns, where, typically, each column refers to a distinct observable or feature. Provided that the data is well formed, such as having zero mean, bound variance, and negligible autocorrelation along each column, the row-space of~$P$ forms the \emph{constituent basis}~$\Pi \in \realsspace^{n\times n}$ with basis vectors~$\pi \in \realsspace^n$ such that \ebox{\Pi = \left(\pi_1, \pi_2, \ldots, \pi_n\right)} (see figure~\ref{fig: r3 point clouds}). If we label the columns in~$P$ by \ebox{\left(p_1, \ldots, p_n\right)}, then the basis~$\Pi$ when materialized onto the constituent axes is simply
\begin{equation*}
  \Pi \; = \;
  \begin{blockarray}{ccccc}
    p_1 & p_2 & \cdots & p_n \\
      \begin{block}{(cccc)c}
      1  &   &        &   & p_1 \\
         & 1 &        &   & p_2 \\
         &   & \ddots &   & \vdots \\
         &   &        & 1 & p_n \\
      \end{block}
  \end{blockarray}
  \quad ,
\end{equation*}
or, \ebox{\Pi = \Identity}. 

Now, when plotted as a scatter plot in~$\realsspace^{n}$, the data in the rows of~$P$ typically form an ellipsoid in the space. When the columns of~$P$ have zero linear correlation, the axes of the ellipsoid align to the constituent axes. As correlation is introduced across the columns of~$P$, the ellipsoid rotates away. Given such a case, we might seek to rotate the ellipsoid back onto the constituent axes. The constituent axes can therefore act as our reference. 

This idea has good flavor, but is incomplete. An axis is different from a vector, for a vector may point along either the positive or negative direction of a particular axis. Therefore, ambiguity remains. To capture this mathematically, let us construct a rotation matrix~$R \in \realsspace^{n\times n}$ such that~$R\transpose$ rotates the data ellipse onto the constituent axes. A modified \SVD will read
\begin{equation*}
  P R = U\: \EigValMatrix^{1/2}\: V\transpose R.
\end{equation*}
If it were the case that \ebox{R\transpose V = \Identity}, then the goal would be achieved. However, \ebox{R\transpose V = \Identity} \emph{is not true}. Instead,
\begin{equation}
  R\transpose V \sim \diagoper\left( \pm 1, \pm 1, \ldots, \pm 1  \right),
  \label{eq: RT V sim diag(signs)}
\end{equation}
is the general case, where the sign on the~$k$th entry depends on the corresponding eigenvector orientation in~$V$. While the following is an oversimplification of the broader concept of orientation, it is instructive to consider that a basis in~$\realsspace^{3}$ may be oriented in a right- or left-handed way. It is not possible to rotate one handedness onto the other because of the embedded reflection. 

It is worth noting that an objection might be raised at this point because the eigenvector matrix~$V$ is itself a unitary matrix with \ebox{\det(V) = +1} and therefore serves as the sought-after rotation. Indeed, \ebox{V\transpose V = \Identity}. The problem is that~$V\transpose V$ is quadratic, so any embedded reflections cancel. In order to identify the vector orientation within~$V$, we must move away from the quadratic form.

%
%

\section{A Representation for Consistent Orientation}
\label{sec: a representation for consistent orientation}

The representation detailed herein uses the constituent basis \ebox{\Pi = \Identity} as the reference basis and seeks to make~(\ref{eq: RT V sim diag(signs)}) an equality for a particular~$V$. To do so, first define the diagonal sign matrix~$S_k \in \realsspace^{n\times n}$ as
\begin{equation}
  S_k = \diagoper\left(1, 1, \ldots, s_k, \ldots, 1 \right),
  \quad\text{where}\quad
  s_k = \pm 1.
  \label{eq: S_k def}
\end{equation} 
The product \ebox{S = \prod S_k} is \ebox{S = \diagoper\left( s_, s_2, \ldots, s_n\right)}, and clearly \ebox{S^2 = \Identity}. Given an entry \ebox{s_k = -1}, the $V S_k$ product reflects the $k$th eigenvector such that \ebox{v_k \rightarrow -v_k} while leaving the other eigenvectors unaltered.  

An equality can now be written with the representation
\begin{equation}
  R\transpose V S = \Identity,
  \label{eq: RT V S = I representation def}
\end{equation}
given suitable rotation and reflection matrices~$R$ and~$S$. %
This representation says that a rotation matrix~$R$ can be found to rotate~$V S$ onto the constituent basis~$\Pi$, provided that a diagonal reflection matrix~$S$ is found to selectively flip the sign of eigenvectors in~$V$ to ensure complete alignment in the end. Once~(\ref{eq: RT V S = I representation def}) is satisfied, the oriented eigenvector basis~$\orientedeigvec$ can be calculated in two ways:
\begin{equation}
  \orientedeigvec = V S = R. 
  \label{eq: V-orient = V S = R}
\end{equation}

To achieve this goal, the first step is to order the eigenvalues in descending order, and reorder the eigenvectors in~$V$ accordingly. (Since the representation~(\ref{eq: RT V S = I representation def}) holds for real symmetric matrices \ebox{M=M\transpose} and not solely for positive semidefinite systems, the absolute value of the eigenvalues should be ordered in descending order.) Indexing into~$V$ below expects that this ordering has occurred. 

The representation~(\ref{eq: RT V S = I representation def}) is then expanded so that each eigenvector is inspected for a reflection, and a rotation matrix is constructed for its alignment to the constituent. The expansion reads
\begin{equation}
  R_n\transpose \ldots R_2\transpose R_1\transpose\: V S_1 S_2 \ldots S_n = \Identity . 
  \label{eq: representation expanded}
\end{equation}
There is a good separation of concerns here: reflections are imparted by the~$S_k$ operators to the right of~$V$, while rotations are imparted by the~$R_k\transpose$ operators to the left. Reflections \emph{do not} preserve orientation whereas rotations \emph{do} preserve orientation. 

The algorithm created to attain a consistent orientation follows the pattern
\begin{equation}
  V \rightarrow V S_1 \rightarrow R_1\transpose V S_1 
                      \rightarrow R_1\transpose V S_1 S_2 
                      \rightarrow \cdots \;.
  \label{eq: program to attain R V S = I}
\end{equation}
For each eigenvector, a possible reflection is computed first, followed by a rotation. The result of the rotation is to align the target eigenvector with its associated constituent basis vector. The first rotation is performed on the entire space~$\realsspace^n$, the second  is performed on the subspace~$\realsspace^{n-1}$ so that the alignment \ebox{v_1\transpose\pi_1 = 1} is preserved, the third rotation operates on~$\realsspace^{n-2}$ to preserve the preceding alignments, and so on. The solution to~(\ref{eq: RT V S = I representation def}) will include~$n$ reflection elements~$s_k$ and \ebox{n(n-1)/2} rotation angles embedded in the~$n$ rotation matrices~$R_k$. 

Before working through the solution, it will be helpful to walk through an example of the attainment algorithm in~$\realsspace^3$.

%
%

\section{\texorpdfstring{Algorithm Example in~$\realsspace^3$}{Algorithm Example in R3}}
\label{sec: algorithm example in r3}

This example illustrates a sequence of reflections and rotations in~$\realsspace^3$ to attain alignment between~$V$ and~$\Pi$, and along the way shows cases in which a reflection or rotation is simply an identity operator. 

\begin{figure*}[t!]\sidecaption
  \centering
  \includegraphics[width=126mm]{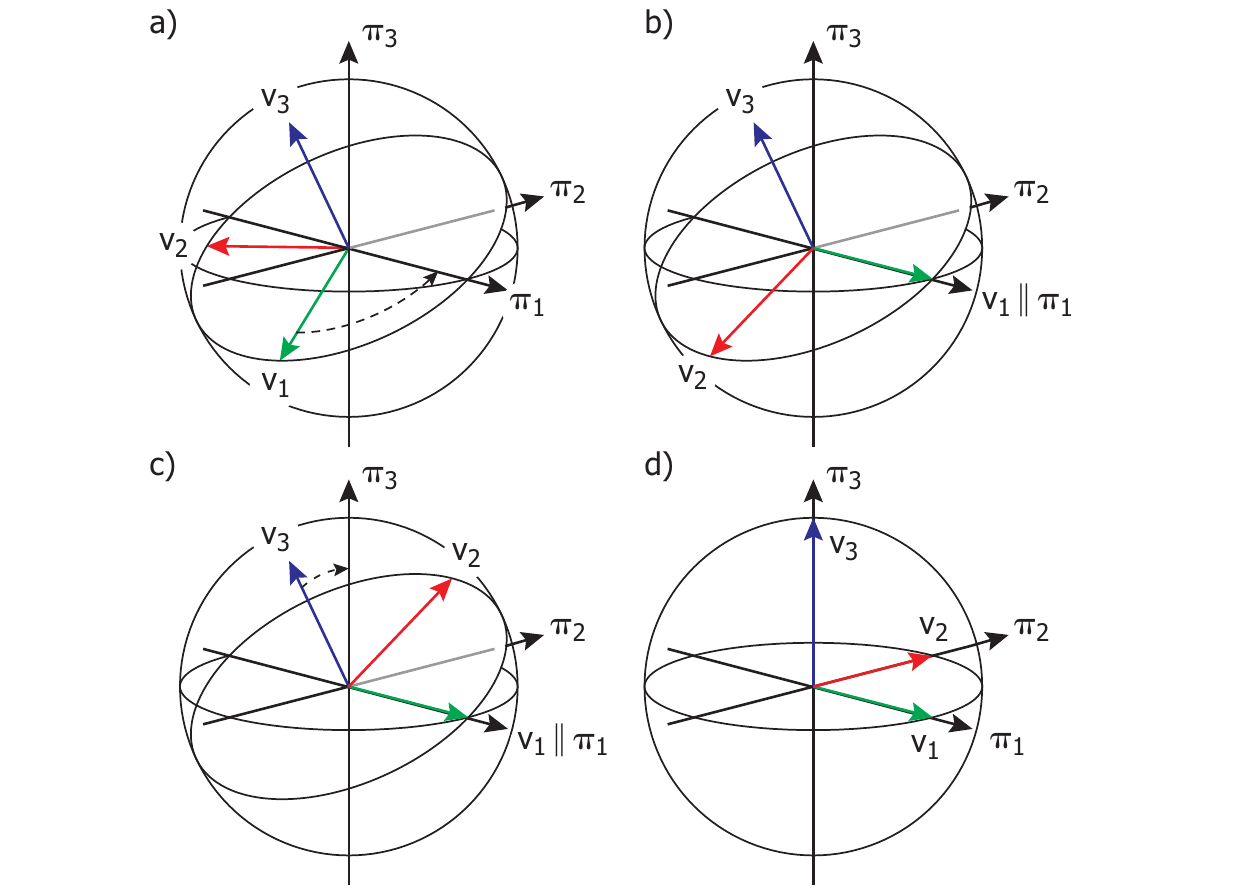}
  \caption{Example of the orientation of a left-handed eigenvector basis~$V$ to the constituent basis~$\Pi$. a) Orthonormal eigenvectors \ebox{(v_1, v_2, v_3)} as they are oriented with respect to the constituent basis vectors \ebox{(\pi_1, \pi_2, \pi_3)}. b) Rotation of~$V$ so that \ebox{v_1\transpose \pi_1 = 1}. c) Reflection of~$v_2$ so that \ebox{(R_1\transpose v_2)\transpose \pi_2 \geq 0}. d) Rotation of $R_1\transpose V S_1$ by~$R_2$ to align~$V$ with~$\Pi$. 
  }
  \label{fig: algorithm example in r3}
\end{figure*}

A example eigenvector~$V$ is created with a left-hand basis, such that an alignment to the constituent basis~$\Pi$ would yield \ebox{R\transpose V = \diagoper(1, -1, 1)}. The precondition of ordered eigenvalues \ebox{\lambda_1 > \lambda_2 > \lambda_3} has already been met. Figure~\ref{fig: algorithm example in r3}~(a) shows the relative orientation of~$V$ to~$\Pi$. 

%

The first step is to consider whether or not to reflect~$v_1$, the purpose being to bring~$v_1$ into the nonnegative side of the halfspace created by the hyperplane perpendicular to~$\pi_1$.  To do so, the relative orientation of \ebox{(v_1, \pi_1)} is measured with an inner product. The reflection entry~$s_1$ is then determined from
\begin{equation*}
  s_1 = \sign\left( v_1\transpose \pi_1\right)
\end{equation*}
(with \ebox{\sign(0) = 1}). From this, $V S_1$ is formed. A consequence of the resultant orientation is that the angle subtended between~$v_1$ and~$\pi_1$ is \ebox{\leq |\pi/2|}. 

The next step is to align~$v_1$ with~$\pi_1$ through a rotation~$R_1$: the result is shown in figure~\ref{fig: algorithm example in r3}~(b). The principal result of \ebox{R_1\transpose V S_1} is that \ebox{v_1\transpose\pi_1 = 1}: the basis vectors are aligned. Significantly, the remaining vectors \ebox{(v_2, v_3)} are perpendicular to~$\pi_1$, and therefore the subsequent rotation needs only to act in the \ebox{(\pi_2, \pi_3) \in \realsspace^2} subspace. 

With this sequence of operations in the full space~$\realsspace^3$ now complete, the subspace~$\realsspace^2$ is treated next. As before, the first step is to inspect the orientation of~$v_2$ with~$\pi_2$. However, care must be taken because~$v_2$ is no longer in its original orientation; instead, it was rotated with the whole basis during $R_1\transpose V$, as shown in figure~\ref{fig: algorithm example in r3}~(b). It is the \emph{new} orientation of~$v_2$ that needs inspection, and that orientation is $R_1\transpose v_2$. As seen in figure~(b), $R_1\transpose v_2$ and~$\pi_2$ point in opposite directions, so eigenvector~$v_2$ needs to be reflected. We have
\begin{equation*}
  s_2 = \sign\left( (R_1\transpose v_2)\transpose \pi_2 \right),
\end{equation*}
which in this example is \ebox{s_2 = -1}. The updated representation of the eigenbasis is \ebox{R_1\transpose V S_1 S_2}, as seen in figure~\ref{fig: algorithm example in r3}~(c). 

To complete work in~$\realsspace^2$, rotation~$R_2$ is found to bring~$R_1\transpose v_2 s_2$ into alignment with~$\pi_2$. The sequence \ebox{R_2\transpose R_1\transpose V S_1 S_2} is illustrated in figure~\ref{fig: algorithm example in r3}~(d). 

The solution needs to work through the last subspace~$\realsspace$ to align~$v_3$. The reflection entry~$s_3$ is evaluated from
\begin{equation*}
  s_3 = \sign\left( (R_2\transpose R_1\transpose v_3)\transpose \pi_3 \right),
\end{equation*}
where \ebox{s_3 = +1} in this example. From this we arrive at \ebox{R_2\transpose R_1\transpose V S_1 S_2 S_3}. Canonically, a final rotation~$R_3$ will carry $R_2\transpose R_1\transpose v_3 s_3$ into~$\pi_3$, but since we are in the last subspace of~$\realsspace$, the operation is a tautology. Still, $R_3$ is included in the representation~(\ref{eq: representation expanded}) for symmetry. 

Lastly, the final, oriented eigenvector basis can be calculated two ways,
\begin{equation}
  \orientedeigvec = V S_1 S_2 S_3 =   R_1 R_2 R_3;
  \label{eq: Vor = V S1 S2 S3 = R1 R2 R3}
\end{equation}
compare with~(\ref{eq: V-orient = V S = R}).

In summary, the algorithm walks down from the full space~$\realsspace^3$ into two subspaces in order to perform rotations that do not upset previous alignments. As this example shows, \ebox{S_1 = S_3 = \Identity}, so there are in fact two reflection operators that are simply identities. Moreover, \ebox{R_3 = \Identity}, since the prior rotations automatically align that last basis vector to within a reflection.

%
%

\section{Algorithm Details}
\label{sec: algorithm details}

The algorithm that implements~(\ref{eq: representation expanded}) is now generalized to~$\realsspace^n$. There are three building blocks to the algorithm: sorting of the eigenvalues in descending order, calculation of the reflection entries, and construction of the rotation matrices. Of the three, construction of the rotation matrices requires the most work and has not been previously addressed in detail.

\subsection{Sorting of Eigenvalues}
\label{subsec: sorting of eigenvalues}

The purpose of sorting the eigenvalues in descending order and rearranging the corresponding eigenvectors accordingly is to ensure that the first rotation operates on the axis of maximum variance, the second rotation on the axis of next-largest variance, and so on, the reason being that the numerical certainty of the principal ellipsoid axis is higher than that of any other axis. Since the rotations are concatenated in this method, an outsize error introduced by an early rotation persists throughout the calculation. By sorting the eigenvalues and vectors first, any accumulated errors will be minimized. 

Since the algorithm works just as well for real-valued symmetric matrices, it is reasonable to presort the eigenvalues in descending order of their absolute value.  From a data-science perspective, we generally expect positive semidefinite systems; therefore a system with negative eigenvalues is unlikely, and if one is encountered, there may be larger upstream issues to address.

\subsection{Reflection Calculation}
\label{subsec: reflection calculation}

Reflection entries are calculated by inspecting the inner product between a rotated eigenvector and its associated constituent basis (several examples are given in the previous section). The formal reflection-value expression for the~$k$th basis vector is
\begin{equation}
  s_k = \sign
    \left( \left( \mathcal{R}_k\transpose v_k \right)\transpose \pi_k \right)
  \label{eq: s_k expr def}
\end{equation}
where
\begin{equation}
  \mathcal{R}_k\transpose \equiv R_{k-1}\transpose \cdots R_2\transpose R_1\transpose \;.
  \label{eq: R_k cal def}  
\end{equation}
Despite of the apparent complexity of these expressions, the practical implementation is simply to inspect the sign of the \ebox{(k, k)} entry in the \ebox{R_{k-1}\transpose\cdots R_1\transpose V} matrix.

\subsection{Rotation Matrix Construction}
\label{subsec: rotation matrix construction}

One way to state the purpose of a particular rotation is that the rotation reduces the dimension of the subspace by one. For an eigenvector matrix~$V$ in~$\realsspace^n$, denoted by~$V_{(n)}$, the rotation acts such that
\begin{equation}
  R_1\transpose \left( V_{(n)} S_1 \right) = 
  \left(
    \begin{array}{cccc}
      1 & & & \\
        & \multicolumn{3}{c}{
            \multirow{3}{36pt}{ $V_{(n-1)}$ }} \\ \\ \\
    \end{array}
  \right) \;.
  \label{eq: R1 V_(n) -> (1, V_(n-1))}
\end{equation}
Here, $S_1$ ensures that the~$(1,1)$ entry on the right is~$+1$. Next, the sign of the~$(2,2)$ matrix entry is inspected and~$s_2$ set accordingly. The next rotation decrements the dimension again such that
\begin{equation}
  R_2\transpose \left( R_1\transpose V_{(n)} S_1 S_2 \right) = 
  \left(
    \begin{array}{cccc}
      1 &   & & \\
        & 1 & & \\
        &   & \multicolumn{2}{c}{ V_{(n-2)} } \\  
    \end{array}
  \right) \;.
  \label{eq: R2 V_(n-1) -> (1, 1, V_(n-2))}
\end{equation}
As before, $S_2$ ensures that the~$(2,2)$ entry on the right is~$+1$. The~$1$ entries along the diagonal reflect the fact that a eigenvector has been aligned with its corresponding constituent basis vector. Subsequent rotations preserve the alignment by operating in an orthogonal subspace.

\subsubsection{Alignment to One Constituent Basis Vector}
\label{subsubsec: alignment to one constituent basis vector}

To simplify the notation for the analysis below, let us collapse~$V$ and all operations prior to the~$k$th rotation into a single working matrix, denoted by~$W_k$. The left-hand side of the two equations above then deals with \ebox{R_1\transpose W_1} and \ebox{R_2\transpose W_2}. The~$\ell$th column of the $k$th working matrix is~$w_{k, \ell}$, and the~$m$th column entry is~$w_{k, \ell, m}$. 

With this notation, the actions of~$R_1$ and~$R_2$ above on the principal working columns~$w_{1,1}$ and~$w_{2,2}$ are
\begin{equation}
  R_1\transpose 
    \left(
      \begin{array}{c}
        w_{1,1,1} \\
        w_{1,1,2} \\
        \vdots \\
        w_{1,1,n}
      \end{array}
    \right) = 
    \left(
      \begin{array}{c}
        1 \\
        0 \\
        \vdots \\
        0
      \end{array}
    \right)
  \;\;\text{and}\;\;
  R_2\transpose 
    \left(
      \begin{array}{c}
        0 \\
        w_{2,2,2} \\
        \vdots \\
        w_{2,2,n}
      \end{array}
    \right) = 
    \left(
      \begin{array}{c}
        0 \\
        1 \\
        \vdots \\
        0
      \end{array}
    \right) \;.
  \label{eq: rotate w into column of zeros with 1 on the pivot}
\end{equation}
Strictly speaking, the~$1$ entries in the column vectors on the right-hand sides are \ebox{w_{k,k}\transpose w_{k,k}}, but since~$V$ is expected to be orthonormal, the inner product is unity. 

Equations like~(\ref{eq: rotate w into column of zeros with 1 on the pivot}) are well known in the linear algebra literature, for instance \cite{Strang:1988,Golub:2014}, and form the basis for the implementation of QR decomposition (see also \cite{NRC:1992}). There are two approaches to the solution: use of Householder reflection matrices, or use of Givens rotation matrices. The Householder approach is reviewed in section~\ref{subsec: householder reflections} below, but in the end it is not suitable in this context. The Givens approach works well, and the solution here has an advantage over the classical formulation because we can rely on \ebox{w_{k,k}\transpose w_{k,k} = 1}.

To simplify the explanation once again, let us focus on the four-dimensional space~$\realsspace^4$. In this case, $R_1\transpose w_{1,1}$ needs to zero out three entries below the pivot, and likewise $R_2\transpose w_{2,2}$ needs to zero out two entries below its pivot. Since one Givens rotation matrix can zero out one column entry at most, a cascade of Givens matrices is required to realize the pattern in~(\ref{eq: rotate w into column of zeros with 1 on the pivot}) and the length of the cascade depends on the dimension of the current subspace. 

A Givens rotation imparts a rotation within an~$\realsspace^2$ plane embedded in a higher-dimensional space. While a simple two-dimensional counterclockwise rotation matrix is
\begin{equation*}
  R(\theta) = 
  \left(
    \begin{array}{cc}
      \cos\theta  &  -\sin\theta \\
      \sin\theta  &   \cos\theta
    \end{array}
  \right),
\end{equation*}
an example Givens rotation matrix in~$\realsspace^4$ is
\begin{equation*}
  R_{\cdot, 3}\left(\theta_{\cdot, 3}\right) = 
  \left(
    \begin{array}{cccc}
      c_3  &     & -s_3  &  \\
           & 1   &       &  \\
      s_3  &     & c_3   &  \\
           &     &       & 1
    \end{array}
  \right)
  \longrightarrow
  \left(
    \begin{array}{cccc}
      \centerdot &            & \centerdot & \\    
                 & \centerdot &            & \\    
      \centerdot &            & \centerdot & \\    
                 &            &            & \centerdot
    \end{array}  
  \right),
%
\end{equation*}
where the dot notation on the right represents the pattern of nonzero entries. The full rotation~$R_1$ is then decomposed into a cascade of Givens rotations:
\begin{equation}
  R_1(\theta_{1,2}, \theta_{1,3}, \theta_{1,4}) = 
    R_{1,2}(\theta_{1,2})\: R_{1,3}(\theta_{1,3})\: R_{1,4}(\theta_{1,4}) \;.
  \label{eq: Givens rotations cascaded into R}
\end{equation}
The order of the component rotations is arbitrary and is selected for best analytic advantage. That said, a different rotation order yields different Givens rotation angles, so there is no uniqueness to the angles calculated below. The only material concern is that the rotations must be applied in a consistent order. 

A advantageous cascade order of Givens rotations to expand the left-hand equation in~(\ref{eq: rotate w into column of zeros with 1 on the pivot}) follows the pattern
\begin{equation}
\underbrace{
    \left(
      \begin{array}{cccc}
        \centerdot & \centerdot &            & \\
        \centerdot & \centerdot &            & \\
                   &            & \centerdot & \\    
                   &            &            & \centerdot
      \end{array}
    \right)
}_{R_{1,2}}
\;
\underbrace{
    \left(
      \begin{array}{cccc}
        \centerdot &            & \centerdot & \\
                   & \centerdot &            & \\    
        \centerdot &            & \centerdot & \\
                   &            &            & \centerdot
      \end{array}
    \right)
}_{R_{1,3}}
\;
\underbrace{
    \left(
      \begin{array}{cccc}
        \centerdot &            &            & \centerdot \\
                   & \centerdot &            & \\
                   &            & \centerdot & \\    
        \centerdot &            &            & \centerdot 
      \end{array}
    \right)
}_{R_{1,4}}
\;
\underbrace{
    \left(
      \begin{array}{c}
        \centerdot \\
                   \\
                   \\
                   \\
      \end{array}
    \right)
}_{\kronvector_{1}}
=
\underbrace{
    \left(
      \begin{array}{c}
        \centerdot \\
        \centerdot \\
        \centerdot \\
        \centerdot
      \end{array}
    \right)
}_{w_{1,1}}
  \;,
  \label{eq: Givens cascade in R4, dot notation}
\end{equation}
where the rotation~$R_1$ has been moved to the other side of the original equation, and~$\kronvector_k$ denotes a vector in~$\realsspace^{n\times 1}$, in which all entries are zero except for a unit entry in the~$k$th row. Each Givens matrix moves some of the weight from the first row of~$\kronvector_1$ into a specific row below, and otherwise the matrices are not coupled. 

In the same vein, Givens rotations to represent~$R_2$ follow the pattern
\begin{equation*}
\underbrace{
    \left(
      \begin{array}{cccc}
        \centerdot &            &            & \\
                   & \centerdot & \centerdot & \\
                   & \centerdot & \centerdot & \\
                   &            &            & \centerdot
      \end{array}
    \right)
}_{R_{2,3}}
\;
\underbrace{
    \left(
      \begin{array}{cccc}
        \centerdot &            &            &            \\
                   & \centerdot &            & \centerdot \\
                   &            & \centerdot &            \\    
                   & \centerdot &            & \centerdot 
      \end{array}
    \right)
}_{R_{2,4}}
\;
\underbrace{
    \left(
      \begin{array}{c}
                   \\
        \centerdot \\
                   \\
                   \\
      \end{array}
    \right)
}_{\kronvector_2}
=
\underbrace{
    \left(
      \begin{array}{c}
                   \\
        \centerdot \\
        \centerdot \\
        \centerdot
      \end{array}
    \right)
}_{w_{2,2}}
  \;.
\end{equation*}
Similar to before, weight in the~$(2,1)$ entry of~$\kronvector_2$ is shifted into the third and fourth rows of the final vector. It is evident that~$R_3$ requires only one Givens rotation, and that for~$R_4$ is simply the identity matrix.

\subsubsection{Solution for Givens Rotation Angles}
\label{subsubsec: solution for givens rotation angles}

Now, to solve for the Givens angles, the cascades are multiplied through. For the~$R_1$ cascade in~(\ref{eq: Givens cascade in R4, dot notation}), multiplying through yields
\begin{equation}
\left(
  \begin{array}{c}
    c_2\: c_3 \: c_4 \\
    s_2\: c_3 \: c_4 \\
    s_3\: c_4        \\
    s_4
  \end{array}
\right)
=
\left(
  \begin{array}{c}
    a_1 \\
    a_2 \\
    a_3 \\
    a_4
  \end{array}
\right),
  \label{eq: trig vector equation def}
\end{equation}
where entries~$a_k$ denote row entries in~$w_{1,1}$ and are used only to further simplify the notation. This is the central equation, and it has a straightforward solution. Yet first, there are several important properties to note: 
\begin{enumerate}

  \item The~$L_2$ norms of the two column vectors are both unity. 
  
  \item While there are four equations, the three angles~$(\theta_2,$ $\theta_3, \theta_4)$ are the only unknowns. The fourth equation is constrained by the~$L_2$ norm.
    
  \item The~$a_1$ entry (or specifically the $w_{1,1,1}$ entry) is nonnegative: \ebox{a_1 \geq 0}. This holds by construction because the associated reflection matrix, applied previously, ensures the nonnegative value of this leading entry. 
  
  \item In order to uniquely take the arcsine, the rotation angles are restricted to the domain \ebox{\theta \in [-\pi/2, \pi/2]}. As a consequence, the sign of each row is solely determined by the sine-function values, the cosine-function values being nonnegative on this domain: \ebox{c_k \geq 0}. This leads to a global constraint on the solution since \ebox{c_2 c_3 c_4 \geq 0} could otherwise hold for pairs of negatively signed cosine-function values. 
  
\end{enumerate}

With these properties in mind, a solution to~(\ref{eq: trig vector equation def}) uses the arcsine method, which starts from the bottom row and works upward. The sequence in this example reads
\begin{equation}
  \begin{aligned}
    \theta_4 & = \sin\inverse \left( a_4 \right)   \\
    \theta_3 & = \sin\inverse \left( a_3 / \cos\theta_4 \right)  \\
    \theta_2 & = \sin\inverse \left( a_2 / \left( \cos\theta_4 \cos\theta_3 \right)\right) ,
  \end{aligned}
  \label{eq: arcsine solution to givens angles}
\end{equation}
where \ebox{\theta_k \in [-\pi/2, \pi/2]}. The sign of each angle is solely governed by the sign of the corresponding~$a_k$ entry. Moreover, other than the first angle, the Givens angles are coupled in the sense that no one angle can be calculated from the entries of working matrix~$W_k$ alone, and certainly not prior to the rotation of~$W_{k-1}$ into~$W_k$. This shows once again that there is no shortcut to working down from the full space through each subspace until~$\realsspace^1$ is resolved. 

Past the first equation, the arguments to the arcsine functions are quotients, and it will be reassuring to verify that these quotients never exceed unity in absolute value. For~$\theta_3$, we can write
\begin{equation*}
  a_3^2 = 1 - a_4^2 - a_2^2 - a_1^2 = \cos^2\theta_4 - (a_1^2 + a_2^2),
\end{equation*}
thus \ebox{a_3^2 \leq \cos^2\theta_4}, so the arcsine function has a defined value. Similarly, for~$\theta_2$ we have
\begin{equation*}
  a_2^2 = 1 - a_4^2 - a_3^2 - a_1^2 = \cos^2\theta_4 - \sin^2\theta_3\cos^2\theta_4 - a_1^2,
\end{equation*}
and therefore \ebox{a_2^2 \leq \cos^2\theta_3\cos^2\theta_4}. Again, the arcsine function has a defined value. This pattern holds in any dimension. 

The preceding analysis is carried out for each eigenvector in~$V$. For each subspace~$k$ there are~$k-1$ angles to resolve. For convenience, the Givens angles can be organized into the upper right triangle of a square matrix. For instance, in~$\realsspace^4$,
\begin{equation}
  \boldsymbol{\theta} \longleftarrow 
  \left(
    \begin{array}{cccc}
      0 & \theta_{1,2} & \theta_{1,3} & \theta_{1,4} \\
        & 0            & \theta_{2,3} & \theta_{2,4} \\
        &              & 0            & \theta_{3,4} \\
        &              &              & 0
    \end{array}
  \right) \;.
  \label{eq: theta mtx storage}
\end{equation}

Lastly, it is worth noting that Dash constructed a storage matrix of angles similar to~(\ref{eq: theta mtx storage}) when calculating the embedded angles in a correlation matrix \cite{Dash:Polar:2004}. His interest, shared here, is to depart from the Cartesian form of a matrix---whose entries in his case are correlation values and in the present case are eigenvector entries---and cast them into polar form. Polar form captures the natural representation of the system.

\subsubsection{Summary}
\label{subsubsec: summary}

This section has shown that each eigenvector rotation matrix~$R_k$ is constructed from a concatenation of elemental rotation matrices called Givens rotations. Each Givens rotation has an angle, and provided that the rotation order is always preserved, the angles are a meaningful way to represent the orientation of the eigenbasis within the constituent basis. 

The representation for~$\boldsymbol{\theta}$ in~(\ref{eq: theta mtx storage}) highlights that there are \ebox{n (n-1) / 2} angles to calculate to complete the full basis rotation in~$\realsspace^n$. It will be helpful to connect the angles to a rotation matrix, so let us define a generator function~$\mathcal{G}(\cdot)$ such that
\begin{equation}
  R = \mathcal{G}\left(\boldsymbol{\theta}\right)
  \quad\text{and}\quad
  R_k = \mathcal{G}\left(\boldsymbol{\theta}, k\right).
  \label{eq: R = G(theta) def}
\end{equation}
The job of the generator is to concatenate Givens rotations in the correct order, using the recorded angles, to correctly produce a rotation matrix, either a matrix~$R_k$ for one vector, or the matrix~$R$ for the entire basis.

%
%

\section{Consideration of Alternative Approaches}
\label{sec: consideration of alternative approaches}

The preceding analysis passes over two choices made during the development of the algorithm. These alternatives are detailed in this section, together with the reasoning for not selecting them.

\subsection{Householder Reflections}
\label{subsec: householder reflections}

Equations nearly identical to~(\ref{eq: rotate w into column of zeros with 1 on the pivot}) appear in well-known linear algebra texts, such as \cite{Golub:2014,Strang:1988,Strang:1986}, typically in the context of QR decomposition, along with the advice that Householder reflections are preferred over Givens rotations because all entries below a pivot can be zeroed out in one step. This advice is accurate in the context of QR decomposition, but does not hold in the current context. 

Recall that in order to align eigenvector~$v_1$ to constituent basis~$\pi_1$, $n-1$ Givens rotations are necessary for a space of dimension~$n$. Arranged appropriately, each rotation has the effect of introducing a zero in the column vector below the pivot. A Householder reflection aligns~$v_1$ to~$\pi_1$ in a single operation by \emph{reflecting}~$v_1$ onto~$\pi_1$ about an appropriately aligned hyperplane in~$\realsspace^n$. In general this hyperplane is neither parallel nor perpendicular to~$v_1$, and therefore neither are any of the other eigenvectors in~$V$. 

The Householder operator is written as
\begin{equation*}
  H \equiv \Identity - 2 u u\transpose,
\end{equation*}
where~$u$ is the Householder vector that lies perpendicular to the reflecting hyperplane. The operator is Hermitian, unitary, involutory, and has \ebox{\det(H) = -1}. In fact, all eigenvalues are~$+1$, except for one that is~$-1$. 

Let us reconsider the eigenvector representation~(\ref{eq: RT V S = I representation def}) and apply a general unitary transform~$U$ to~$V$, as in
\begin{equation}
  U\transpose V S = \Identity.
  \label{eq: UT V S = I alternative}
\end{equation}
(Unitary transform~$U$ is not to be confused with the matrix of projected data~$U$ in \SVD equation~(\ref{eq: SVD decomp def}).) The equality holds for any unitary operator, so Householder reflections are admitted. Let us then compare expanded representations in the form of~(\ref{eq: representation expanded}) for both rotation and reflection operators to the left of~$V$:
\begin{equation*}
  \begin{aligned}
    \underbrace{R_n\transpose \ldots R_2\transpose R_1\transpose}_{\text{basis rotations}}\;
    \;\;\;V
    \underbrace{S_1 S_2 \ldots S_n}_{\text{eigenvector reflections}} & = \Identity \\[9pt]
    \underbrace{H_n\transpose \ldots H_2\transpose H_1\transpose}_{\text{basis reflections}}\; 
    \;\;\;V 
    \underbrace{S_1 S_2 \ldots S_n}_{\text{eigenvector reflections}} & = \Identity     
  \end{aligned}
\end{equation*}
Even though~(\ref{eq: UT V S = I alternative}) holds for Householder reflections, the equality itself \emph{is not the goal}; rather, the \emph{interpretability} of the eigenvector orientation is paramount. In the rotations method, reflections only reflect one eigenvector at a time whereas rotations transform the basis as a whole. In contrast, the reflections method intermingles eigenvector reflections with basis reflections, thereby disrupting the interpretation of the basis orientation. 

%
\begin{figure*}[t]\sidecaption
  \centering
  \includegraphics[width=126mm]{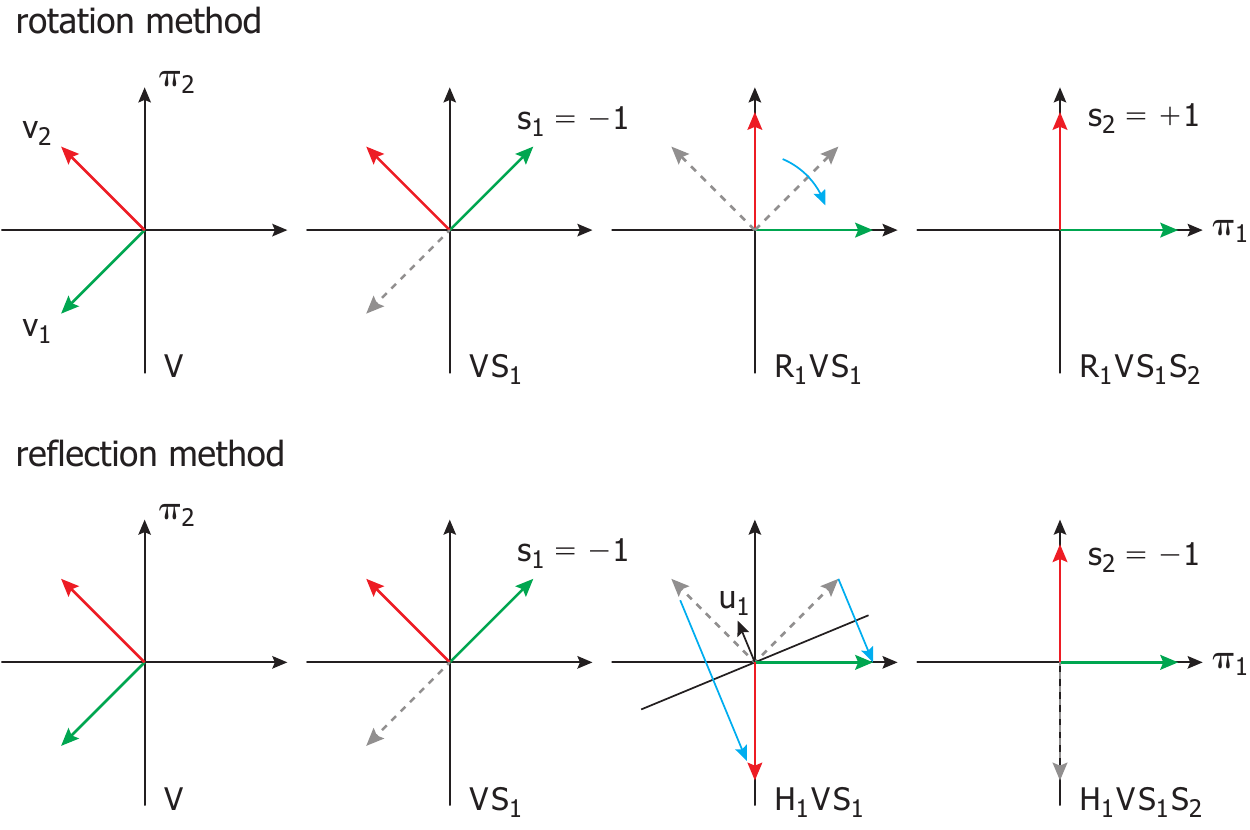}
  \caption{Rotation versus reflection transformation comparison. \emph{Top}, the rotation method produces the sign tuple \ebox{(s_1, s_2) =} \ebox{(-1, +1)}. \emph{Bottom}, the Householder reflection method produces the sign tuple \ebox{(-1, -1)}. 
  }
  \label{fig: rotation vs reflection}
\end{figure*}

A comparison between the rotation and reflection methods is illustrated in figure~\ref{fig: rotation vs reflection}. Here, the eigenvector basis spans~$\realsspace^2$ with an initial orientation of
\begin{equation}
  V =
  \frac{1}{\sqrt{2}}
  \left(
    \begin{array}{cc}
      -1  & -1  \\
      -1  & +1
    \end{array}
  \right) \;.
  \label{eq: V for R2 example}
\end{equation}
Following the upper pane in the figure, where the reflection method is illustrated, it is found that \ebox{v_1\transpose\pi_1 < 0} and therefore \ebox{s_1 = -1} in order to reflect~$v_1$. Once~$V S_1$ is completed, a clockwise rotation of~$45\degrees$ is imparted by~$R_1$ to align~$v_1$ with~$\pi_1$, thus completing $R_1\transpose V S_1$. Lastly, \ebox{v_2\transpose \pi_2 > 0} (trivially so) and therefore \ebox{s_2 = +1}. The tuple of sign reflections for the eigenvectors in~$V$ is therefore \ebox{(s_1, s_2) = (-1, +1)}. 

A parallel sequence is taken along the lower pane, where a Householder reflection is used in place of rotation. Expression~$V S_1$ is as before, but to align~$v_1$ with~$\pi_1$ a reflection is used. To do so, Householder vector~$u_1$ is oriented so that a reflection plane is inclined by~$22.5\degrees$ from the~$\pi_1$ axis. Reflection of~$V S_1$ about this plane indeed aligns~$v_1$ to~$\pi_1$ but has the side effect of reflecting~$v_2$ too. As a consequence, when it is time to inspect the orientation of~$v_2$ with respect to~$\pi_2$, we find that \ebox{s_2 = -1}. Therefore the tuple of eigenvector signs is now \ebox{(-1, -1)}. 

Our focus here is not to discuss which tuple of signs is ``correct,'' but rather on how to interpret the representation. The choice made in this article is to use rotations instead of reflections for the basis transformations so that the basis orientation is preserved through this essential step. The use of Householder reflections, by contrast, scatters the basis orientation after each application.

\subsection{Arctan Calculation Instead of Arcsine}
\label{subsec: arctan calculation instead of arcsine}

A solution to~(\ref{eq: trig vector equation def}) is stated above by equations~(\ref{eq: arcsine solution to givens angles}) in terms of arcsine functions, and a requirement of this approach is that the range of angles be restricted to \ebox{\theta_k \in [-\pi/2, \pi/2]}. While this is the preferred solution, there is another approach that uses arctan functions instead of arcsine. The arctan solution to~(\ref{eq: trig vector equation def}) starts at the top of the column vector and walks downward, the angles being
\begin{equation}
  \begin{aligned}  
    \theta_2 & = \arctantwo\left( a_2, a_1 \csc \theta_1 \right) \\
    \theta_3 & = \arctantwo\left( a_3, a_2 \csc \theta_2 \right) \\
    \theta_4 & = \arctantwo\left( a_4, a_3 \csc \theta_3 \right) \;.
  \end{aligned}
  \label{eq: arctan solution to givens angles}
\end{equation}
The first angle is defined as \ebox{\theta_1 \equiv \pi/2} (which is done solely for equation symmetry), and the~$\arctantwo$ function is the four-quadrant form \ebox{\arctantwo(y, x)}, in which the admissible angular domain is \ebox{\theta \in [-\pi, \pi]}. It would appear that the arctan method is a better choice. 


The difficulty arises with edge cases. Let us consider the vector
\begin{equation*}
  a = \left(1/\sqrt{2}, 1/\sqrt{2}, 0, 0\right)\transpose.
\end{equation*}
The sequence of arctan angles is then
\begin{equation*}
  \begin{aligned}
    \theta_1 & = \pi / 2 \\
    \theta_2 & = \arctantwo\left( 1/\sqrt{2}, 1/\sqrt{2} \times\csc\theta_1 \right) = \pi/4 \\
    \theta_3 & = \arctantwo\left( 0, 1/\sqrt{2} \times\csc\theta_2 \right) = 0 \\
    \theta_4 & = \arctantwo\left( 0, 0 \times \csc\theta_3 \right) \;.
  \end{aligned}
\end{equation*}
The final arctan evaluation is not numerically stable. For computer systems that support signed zero, the arctan can either be \ebox{\arctantwo(0,0) \rightarrow 0} or \ebox{\arctantwo(0, -0) \rightarrow} \ebox{-\pi}, see \cite{numpy:arctan2,ieee:arctan2}. The latter case is catastrophic because product~$V S_1$ has forced the sign of~$v_{1,1}$ to be positive, yet now~$R_1\transpose V S_1$ \emph{flips} the sign since \ebox{\cos\pi = -1}. After~$R_1\transpose V S_1$ is constructed, the first eigenvector is not considered thereafter, so there is no chance, unless the algorithm is changed, to correct this spurious flip. 

Therefore, to avoid edge cases that may give false results, the arcsine method is preferred.

%
%

\section{Application to Regression}
\label{sec: application to regression}

Returning now to the original motivation for a consistently oriented eigenvector basis, the workflow for regression and prediction on the eigenbasis becomes:
\begin{enumerate}
  
  \item \SVD: Find~$(V, \EigValMatrix)$ from in-sample data $P_{\text{in}}$ such that
  \begin{equation*}
    \ebox{P_{\text{in}} = U \EigValMatrix^{1/2}\: V\transpose}.
  \end{equation*}

  \item Orient: Find~$(R, S)$ from~$V$ such that 
  \begin{equation*}
    R\transpose\: V S = \Identity.
  \end{equation*}
  
  \item Regression: Find~$(\hat{\beta}, \epsilon)$ from \ebox{(y, U, S, \EigValMatrix)} such that
  \begin{equation}
    y = U S\: \EigValMatrix^{1/2}\: \hat{\beta} + \epsilon.
    \label{eq: regression with oriented basis}
  \end{equation}  

  \item Prediction: Calculate~$\Expectation y$ from \ebox{(P_{\text{out}}, R, \hat{\beta})}, where $P_{\text{out}}$ is out-of-sample data, with
  \begin{equation}
    \Expectation y = P_{\text{out}}\: R\; \hat{\beta}.
    \label{eq: E y = P R beta} 
  \end{equation}

\end{enumerate}
The form of~(\ref{eq: regression with oriented basis}) here assumes that there is no exogenous correction to the direction of~$y$ in response to a change of the factors along the eigenbasis, as discussed in section~\ref{sec: overview}, \emph{Overview}. 

The workflow highlights that both elements of the orientation solution, rotation~$R$ and reflection~$S$, are used, although for different purposes. The regression (\ref{eq: regression with oriented basis}) requires the reflection matrix~$S$ to ensure that the signs of~$\hat{\beta}$ faithfully align to the response of~$y$. Predictive use with out-of-sample data requires the~$\hat{\beta}$ estimate, as expected, but also the rotation matrix~$R$. The rotation orients the constituent basis, which is observable, into the eigenbasis, which is not. 

Without dimension reduction, the rotation in~(\ref{eq: E y = P R beta}) is associative:
\begin{equation*}
  \bigl( P_{\text{out}} R \bigr) \hat{\beta} = P_{\text{out}} \left( R \hat{\beta} \right)
  \quad\text{for no dimension reduction.}
\end{equation*}
However, associativity breaks with dimension reduction. Principal components analysis, for instance, discards all but the top few eigenvector components (as ranked by their corresponding eigenvalue) and uses the remaining factors in the regression. The number of entries in $\hat{\beta}$ equals the number of remaining components, not the number of constituent components. In this case, only \ebox{\left( P_{\text{out}} R \right)} can be used. Typically, online calculation of this product is simple because the calculation is updated for each new observation. Rather than being an~$\realsspace^{m\times n}$ matrix, out-of-sample updates are $P_{\text{out}} \in \realsspace^{1\times n}$, which means that~$P_{\text{out}}\: R$ requires only a \BLAS level-2 (vector-matrix) call on an optimized system \cite{blas:live}.

%
%

\section{Treatment of Evolving Data}
\label{sec: treatment of evolving data}

\begin{figure*}[t]\sidecaption
  \centering
  \includegraphics[width=126mm]{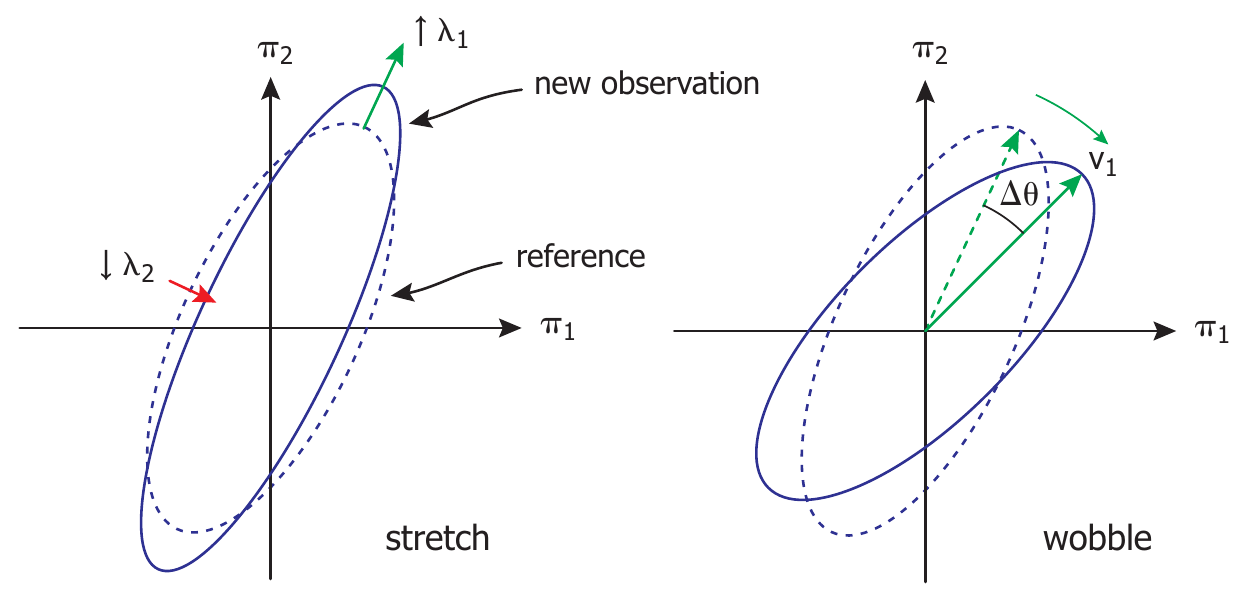}
  \caption{Orthogonal modes of variation of an eigensystem. \emph{Left}, in one variation mode, the \emph{eigenvalues} of a new observation differ from those of the past, thereby stretching and/or compressing the point-cloud ellipse. \emph{Right}, in another variation mode, it is the \emph{eigenvectors} that vary in this way, thereby changing the orientation of the point-cloud ellipse with respect to the constituent basis.
  }
  \label{fig: eigen risks}
\end{figure*}

Eigenvectors and values calculated from data are themselves sample estimates subject to uncertainties based on the particular sample at hand, statistical uncertainties based on the number of independent samples available in each dimension\footnote{See \cite{Meucci:2007}, Section 4.2.2, Dispersion and Hidden Factors.} as well as on assumptions about the underlying distribution,\footnote{See \cite{Meucci:2007}, Section 4.3, Maximum Likelihood Estimators.} and scale-related numerical uncertainties based on the degree of colinearity among factors. All of these uncertainties exist for a single, static dataset.

Additional uncertainty becomes manifest when data evolves because the estimate of the eigensystem will vary at each point in time, this variation being due in part to sample noise, and possibly due to nonstationarity of the underlying random processes. 

In the presence of evolving data, therefore, the eigensystem fluctuates, even when consistently oriented. It is natural, then, to seek statistics for the location and dispersion of the eigensystem. To do so, two orthogonal modes of variation are identified (see figure~\ref{fig: eigen risks}): stretch variation, which is tied to change of the eigenvalues; and wobble variation, which is tied to change of the eigenvectors. Each is treated in turn.

\subsection{Stretch Variation}
\label{subsec: stretch variation}

Stretch variation is nominally simple because eigenvalues~$\lambda_i$ are scalar, real numbers. The mean and variance follow from the usual sample forms of the statistics. For an ensemble of~$N$ samples, the average eigenvalues are simply
\begin{equation}
  \bar{\EigValMatrix} = (1/N) \sum_{i=1}^N \EigValMatrix[i].
  \label{eq: eigenvalue average}
\end{equation}
A challenge for the variance statistic is that eigenvalues are not independent since
\begin{equation*}
  \det\left(V\right) = \lambda_1 \lambda_2 \cdots \lambda_n = +1.
  \label{eq: non-independence of eigenvalues}
\end{equation*}
How the covariance manifests itself is specific to the dataset at hand.

\subsection{Wobble Variation}
\label{subsec: wobble variation}

\begin{figure*}[t]\sidecaption
  \centering
  \includegraphics[width=126mm]{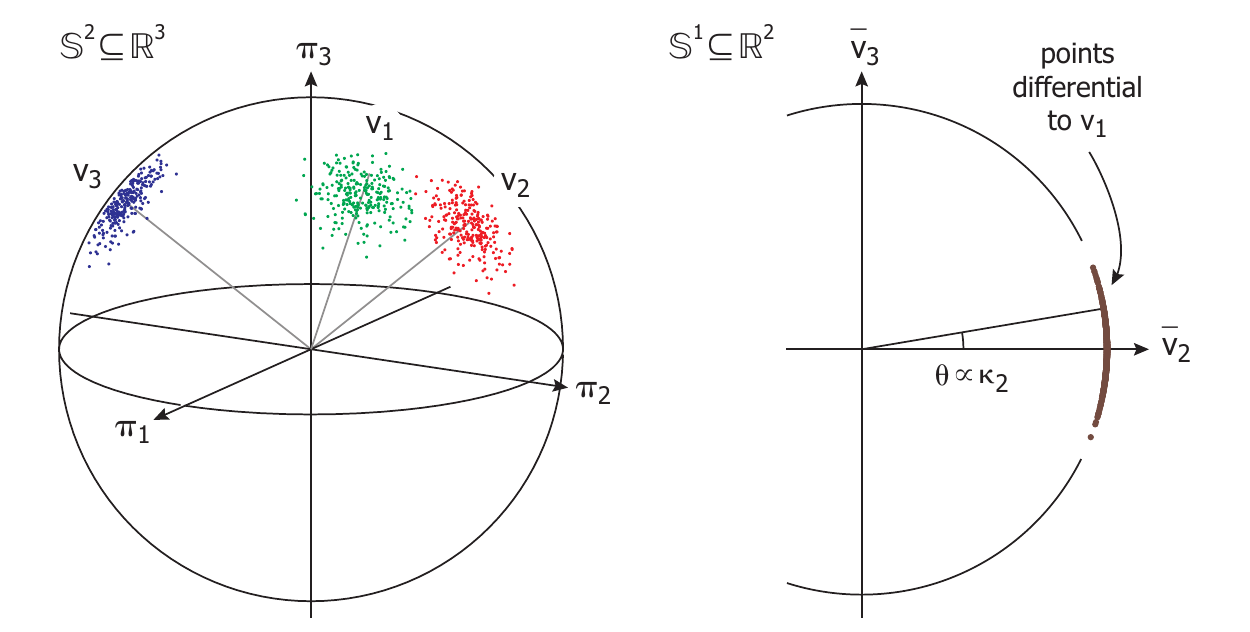}
  \caption{Illustration of eigenvector wobble in~$\realsspace^3$ and its~$\realsspace^2$ subspace. \emph{Left}, points illustrating  eigenvector wobble for a stationary process. The points lie on the hypersphere~$\spherespace^2$, and there is one constellation for each eigenvector. Samples were drawn from a von~Mises--Fisher distribution with concentration \ebox{\kappa = 100}. In this illustration, the point constellations are identical, each having been rotated onto the mean direction of each vector. In general the point constellations are mixtures having different component dispersions. \emph{Right}, eigenvector~$v_2$ may wobble in the \ebox{(\bar{v}_2, \bar{v}_3)} plane independently from~$v_1$ wobble. Points here illustrate wobble in this~$\realsspace^2$ subspace.  
  }
  \label{fig: vmf distributions}
\end{figure*}

Quantification of wobble variation requires a different toolset because eigenvectors are \emph{vectors}, and these orthonormal vectors point onto the surface of a hypersphere such that \ebox{v_k \in \spherespace^{n-1} \subseteq \realsspace^n}. Figure~\ref{fig: vmf distributions}~\emph{left} illustrates a case in point for~$\realsspace^3$: a stationary variation of the eigenvectors forms point constellations on the surface of the sphere~$\spherespace^2$, one constellation for each vector.\footnote{%
Figure~\ref{fig: vmf distributions} was drawn using point constellations drawn from the von~Mises--Fisher distribution \cite{Mardia:2000} with concentration \ebox{\kappa = 100}. Sampling this distribution in~$\realsspace^2$ and~$\realsspace^3$ was done using the \VonMisesFisher class in the \TensorFlow Probability package \cite{tensorflow:probability:2019}, which belongs to the \TensorFlow platform \cite{tensorflow:2015}. The \VonMisesFisher class draws samples using the nonrejection based method detailed by Kurz and Hanebeck \cite{Kurz:2015}. \label{fn: tensorflow}
}
The field of directional statistics informs us regarding how to define mean direction (the vector analogue of location) and dispersion in a consistent manner for point constellations such as these \cite{Mardia:2000,Jammal:2001,Ley:2017}. Directional statistics provides a guide on how to use the vector information available from \ebox{\orientedeigvec} and embedded angle information recorded in~\ebox{\boldsymbol{\theta}}, see~(\ref{eq: theta mtx storage}). Recent work focuses on the application of these statistics to machine learning \cite{Sra:Directional:2018}. Nonetheless, application of directional statistics to eigenvector systems appears to be underrepresented in the literature. 

The following discussion only applies to eigensystems that have been consistent orientated.

%
%

\subsubsection{Estimation of Mean Direction}
\label{subsubsec: estimation of mean direction}

Both mean direction and directional dispersion are measurable statistics for eigensystems, and determining the mean direction is the simpler of the two. The reason this is so is because the eigenvectors are orthogonal in every instance of an eigenbasis, thus the mean directions must also be orthogonal. Consequently, the eigenvectors can be treated uniformly. 

It is a tenet of directional statistics that the \emph{mean direction} is calculated from unit vectors, not from their angles.\footnote{%
Simple examples exist in the cited literature that show that an angle-based measure is not invariant to the choice a ``zero'' angle reference. 
}
For an ensemble of~$N$ unit vectors~$x[i]$, the mean direction \ebox{\bar{x}} is defined as
\begin{equation}
  \bar{x} \;\equiv\; \ltwonorm{x_S} \inverse \; x_S
  \quad\text{where}\quad
  x_S = \sum_{i=1}^N x[i] \;,
  \label{eq: mean direction x-bar def}
\end{equation}
with the caveat that the mean direction is undefined for \ebox{\ltwonorm{x_S} = 0}. The \emph{resultant vector}~$x_S$ is the vector sum of component vectors~$x[i]$ and has length~$\ltwonorm{x_S}$ under an~$L_2$ norm (see figure~\ref{fig: resultant vectors}). The mean direction is thus a unit-vector version of the resultant vector. 

Extending this construction for mean direction to an ensemble of oriented eigenvector matrices~$\orientedeigvec$, the mean location of the eigenvectors is
\begin{equation}
  \bar{\orientedeigvec} = \ltwonorm{\orientedeigvec}\inverse \; \orientedeigvec_S 
  \quad\text{where}\quad
  \orientedeigvec_S = \sum_{i=1}^N \orientedeigvec[i] \;,
  \label{eq: mean direction V-bar def}
\end{equation}
and
\begin{equation}
  \ltwonorm{\orientedeigvec} 
  \equiv 
  \diagoper \Bigl( \ltwonorm{v_{S,1}} , \ltwonorm{v_{S,2}}, \ldots, \ltwonorm{v_{S,N}} \Bigr) 
  \label{eq: mean direction V-bar normalization}
\end{equation}
for normalization. 

Using the methodologies above, the average basis~$\bar{\orientedeigvec}$ can be rotated onto~$\Identity$ with a suitable rotation matrix \ebox{R\transpose \bar{\orientedeigvec} = \Identity}, and in doing so, the Cartesian form of~$\bar{\orientedeigvec}$ can be converted into polar form. Using the generator function~$\mathcal{G}()$ from~(\ref{eq: R = G(theta) def}) to connect the two, we have
\begin{equation}
  \bar\orientedeigvec = \mathcal{G}\left( \bar{\boldsymbol{\theta}} \right).
  \label{eq: polar form of bar-V-oriented}
\end{equation}
Note that \ebox{\bar{\boldsymbol{\theta}}} \emph{is not} the arithmetic average of angles \ebox{\boldsymbol{\theta}[i]} but is exclusively derived from~$\bar{\orientedeigvec}$. In fact, an alternative way to express the basis sum in~(\ref{eq: mean direction V-bar def}) is to write
\begin{equation}
  \orientedeigvec_S = \sum_{i=1}^N \mathcal{G}\Bigl( \boldsymbol{\theta}[i] \Bigr),
  \label{eq: oriented V from generator and theta}
\end{equation}
and from this we see the path of analysis: 
\begin{equation*}
	\boldsymbol{\theta}[i] \rightarrow \orientedeigvec[i] \rightarrow \bar{\orientedeigvec} 			\rightarrow \bar{\boldsymbol{\theta}} .
\end{equation*}

\subsubsection{Estimation of Dispersion}
\label{subsubsec: estimation of dispersion}

Dispersion estimation is more involved for an ensemble of eigenbases than for an ensemble of single vectors because both common and differential modes of variation exist. Referring to figure~\ref{fig: vmf distributions}~\emph{left}, consider a case where the only driver of directional variation for the eigenbasis is a change of the pointing direction of~$v_1$. As~$v_1$ scatters about its mean direction, vectors~$v_{2,3}$ will likewise scatter, together forming three constellations of points, as in the figure. Thus, wobble in~$v_1$ imparts wobble in the other vectors. Since dispersion is a scalar independent of direction, the dispersion estimates along all three directions in this case are identical. 

Yet, there is another possible driver for variation that is orthogonal to~$v_1$, and that is motion within the \ebox{(v_2, v_3)} plane. Such motion is equivalent to a pirouette of the \ebox{(v_2, v_3)} plane about~$v_1$. Taking \ebox{(\bar{v}_2, \bar{v}_3)} as a reference, variation within this~$\realsspace^2$ subspace needs to be estimated, see figure~\ref{fig: vmf distributions}~\emph{right}. When viewed from~$\realsspace^n$, however, the point constellation distribution about~$\bar{v}_2$ is a mixture of variation drivers.

Regardless of the estimation technique for dispersion, it is clear that the pattern developed in section~\ref{subsec: rotation matrix construction}, \emph{Rotation Matrix Construction}, for walking through subspaces of~$\realsspace^n$ to attain \ebox{R\transpose V S = \Identity} occurs here as well. Thus, rather than calculate a dispersion estimate for each eigenvector, since they are mixtures, a dispersion estimate is made for each subspace.

Returning to the generator for subspace~$k$ in~(\ref{eq: R = G(theta) def}), given an ensemble~$\boldsymbol{\theta}[i]$ the ensemble of the~$k$th subspace is
\begin{equation}
  R_k[i] = \mathcal{G}\bigl(\boldsymbol{\theta}[i], k\bigr),
  \label{eq: Rk[i] = G( theta[i], k) }
\end{equation}
where \ebox{R_k \in \realsspace^{n\times n}}. From~$R_k[i]$ the~$k$th column is extracted, and the lower~$k$ entries taken to produce column vector \ebox{x_k[i] \in \realsspace^{(n-k+1)\times 1}}. Figuratively,
\begin{equation}
  R_k = 
  \left(
    \begin{array}{cccc}
      \;1\; & \;0\; & \;0\; & \;0\; \\
      0 & \tikzmarkin{a}(0.20,-0.20)(-0.2,0.3)    
          \centerdot & \centerdot & \centerdot \\
      0 & \centerdot & \centerdot & \centerdot \\
      0 & \centerdot \tikzmarkend{a}
                     & \centerdot & \centerdot 
    \end{array}
  \right)
  \;\longrightarrow\;
  x_k =
  \left(
    \begin{array}{c}
      \tikzmarkin{b}(0.20,-0.20)(0.05,0.3)  
      \;\;\;\centerdot\;\;\; \\
      \centerdot \\
      \centerdot \tikzmarkend{b}
    \end{array}
  \right)  \;.
  \label{eq: Rk to ak}
\end{equation}
The reason behind of removing the top zero entries from the column vector is that parametric models of dispersion are isotropic in the subspace, and therefore keeping dimensions with zero variation would create an unintended distortion. 

\begin{figure*}[t]\sidecaption
  \centering
  \includegraphics[width=126mm]{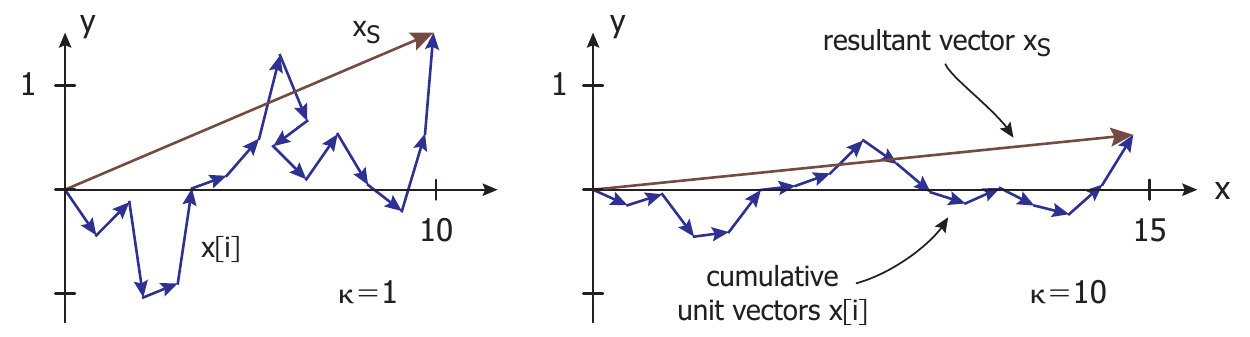}
  \caption{The resultant vector \ebox{x_S \equiv \sum_i x[i]} for two levels of concentration, \ebox{\kappa = (1, 10)}. Lower concentration, which is akin to higher dispersion, leads to a smaller expected norm of the resultant vector \ebox{\ltwonorm{x_S}}, whereas higher concentration leads to a larger expected norm. Sample vectors in~$\realsspace^2$ were generated from the von~Mises--Fisher distribution.  (Note that the difference in horizontal and vertical scales distorts the length of the unit vectors in the plots.)
  }
  \label{fig: resultant vectors}
\end{figure*}

Now that the vectors of interest~$x_k[i]$ have been identified, their dispersions can be considered. Unlike the mean direction, there is no one measure of dispersion. Two approaches are touched upon here, one being model free and the other based on a parametric density function. Both approaches rely on the resultant vector~$x_S$. 

Returning to figure~\ref{fig: resultant vectors}, the two panels differ in that on the left there is a higher degree of randomness in the pointing directions of the unit vectors along the sequence, while on the right there is a lower degree of randomness. The consequence is that the expected \emph{lengths} of the result vectors are different: the lower the randomness, the longer the expected length. Define the \emph{mean resultant length} for~$N$ samples as
\begin{equation}
  \bar{r} = \ltwonorm{x_S} / N, \quad 0 \leq \bar{r} \leq 1.
  \label{eq; mean resultant length def}
\end{equation}
For distributions on a circle, Mardia and Jupp define a model-free \emph{circular variance} as \ebox{\bar{V} \equiv 1 - \bar{r}} (see \cite{Mardia:2000}). Among other things, this variance is higher for more dispersed unit vectors and lower for less dispersed vectors. Yet, after the generalization to higher dimensions, it remains unclear how to compare the variance from one dimension to another. 

As an alternative, a model-based framework starts by positing an underlying parametric probability distribution and seeks to define its properties and validate that they meet the desired criteria. In the present case, eigenvectors are directional, as opposed to axial; their evolution-induced scatter is probably unimodal, at least across short sequences; and their dimensionality is essentially arbitrary. An appropriate choice for an underlying distribution therefore is the von~Mises--Fisher (\vMF) distribution. The \vMF density function in~$\realsspace^n$ is \cite{Sra:Directional:2018}
\begin{equation}
  \begin{gathered}  
	  \vmfpdf(x; n, \vmfmean, \vmfconc) = c_n(\vmfconc)\; e^{\vmfconc \vmfmean\transpose x} 
	  \quad\text{where} \\
	  c_n(\vmfconc) = \frac{ \vmfconc^{n/2-1} }{ (2\pi)^{n/2}\: \modbesselfcxnfirst{n/2-1}(\vmfconc) } \;.
  \end{gathered}
  \label{eq: vmf density def}
\end{equation}
The single argument to the density function is the vector \ebox{x \in \realsspace^{n\times 1}}, and the parameters are the dimension~$n$, mean direction \ebox{\vmfmean \in \realsspace^{n\times 1}}, and concentration \ebox{\vmfconc \in [0, \infty)}. The mean direction parameter is the same as in~(\ref{eq: mean direction x-bar def}): \ebox{\vmfmean = \bar{x}}. The concentration~$\vmfconc$ parameter, a scalar, is like an inverse variance. For \ebox{\vmfconc = 0}, unit vectors~$x$ are uniformly distributed on the hypersphere \ebox{\spherespace^{n-1}}, whereas with \ebox{\vmfconc \rightarrow \infty}, the density concentrates to the mean direction~$\vmfmean$. Maximum likelihood estimation yields
\begin{equation*}
  \hat{\vmfmean} = \ltwonorm{x_S}\inverse\; x_S
  \quad\text{and}\quad
  \hat{\vmfconc} = A_n\inverse\left(\bar{r}\right),
\end{equation*}
where the nonlinear function~$A_n(\cdot)$ is 
\begin{equation*}
  A_n\left(\vmfconc\right) 
    = \frac{ \modbesselfcxnfirst{n/2}(\vmfconc) }{ \modbesselfcxnfirst{n/2-1}(\vmfconc) } 
    = \bar{r}
\end{equation*}
and~$\modbesselfcxnfirst{\alpha}(x)$ is the modified Bessel function of the first kind. A useful approximation is
\begin{equation*}
  \hat{\vmfconc} \simeq \frac{\bar{r} \left(n - \bar{r}^2 \right)}{ 1 - \bar{r}^2 }. 
\end{equation*}
The key aspect here is the connection of the concentration parameter~$\vmfconc$ to the mean resultant length~$\bar{r}$ from~(\ref{eq; mean resultant length def}). Notice as well that~$\vmfconc$ is related to the dimensionality~$n$, indicating that concentration values across different dimensions cannot be directly compared. 

Sampling from the \vMF distribution is naturally desirable. Early approaches used rejection methods (see the introduction in \cite{Kurz:2015}), these being simpler, yet the runtime is not deterministic. Kurz and Hanebeck have since reported a stochastic sampling method that is deterministic, see \cite{Kurz:2015}, and it is this method that is implemented in \TensorFlow (see footnote~\ref{fn: tensorflow} on page~\pageref{fn: tensorflow}). The samples in figure~\ref{fig: vmf distributions} were calculated with \TensorFlow. 

To conclude this section, in order to find the underlying drivers of variation along a sequence of eigenbasis observations, concentration parameters~$\vmfconc$, or at least the mean resultant lengths~$\bar{r}$, should be computed for each descending subspace in~$\realsspace^n$. The~$n$ concentration parameters might then be stored as a vector,
\begin{equation}
  \vmfconc_{\text{basis}} = \left(\vmfconc_1, \vmfconc_2, \ldots, \vmfconc_n  \right) \transpose.
  \label{eq: vmf conc vector over dims}
\end{equation}

\subsection{Rank-Order Change of Eigenvectors}
\label{subsec: rank-order change of eigenvectors}

Along an evolving system, the rank order of eigenvectors may change. This paper proposes no mathematically rigorous solution to dealing such change---it is unclear that one exists, nor is its absence a defect of the current work---yet there are a few remarks that may assist the analyst to determine whether rank-order change is a property of the data itself or a spurious effect that results from an upstream issue. 

Eigenvectors are ``labeled'' by their eigenvalue: in the absence of eigenvalue degeneracy, each vector is mapped one-to-one with a value. The concept of rank-order change of eigenvectors along evolving data is only meaningful if a second label can be attached to the vectors such that the second label does not change when the first does. Intuitively (for evolving data), a second label is the pointing direction of the vectors. This choice is fraught: Two time series are not comparable if they begin at different moments along a sequence because the initial labeling may not match. Additionally, validation of  secondary labeling stability needs to be performed, and revalidation is necessary for each update because stability cannot be guaranteed. 

With an assumption that a second labeling type is suitably stable, rank-order changes may still be observed. Before it can be concluded that the effect is a trait of the data, modeling assumptions must be considered first. There is a significant assumption embedded in the estimation of the eigensystem itself: using an \code{eig(..)} function call presupposes that the underlying statistical distribution is Gaussian. However, data is often heavy tailed, so copula or implicit methods are appropriate since Gaussian-based methods are not robust (see \cite{Meucci:2007}). Another modeling concern is the number of independent samples per dimension used in each panel: too few samples inherently skews the eigenvalue spectrum, and when coupled with sample noise, may induce spurious rank-order changes (again, see \cite{Meucci:2007}). Lastly, \PCA ought only be performed on homogeneous data categories because the mixing of categories may well lead to spurious rank-order changes. Better to apply \PCA independently to each category and then combine. 

If rank-order change is determined to be a true trait of the data, then axial statistics is used in place of directional statistics. Referring to figure~\ref{fig: vmf distributions} above, two or more constellations will have some of their points mirror-imaged through the origin, creating a ``barbell'' shape along a common \emph{axis}. The juxtaposition of the Watson density function, which is an axial distribution, to the von Mises--Fisher distribution, 
\begin{equation*}
	\begin{aligned}	
	    p_{\text{watson}}(x; \vmfmean, \vmfconc) & = z_1\: e^{\vmfconc \left(\vmfmean\transpose x\right)^2}
	    \quad\text{vs}\\
	    \vmfpdf(x; \vmfmean, \vmfconc)           & = z_2\: e^{\vmfconc \left(\vmfmean\transpose x\right)}  \;,
	\end{aligned}
\end{equation*}
shows that the axial direction is squared in order to treat the dual-signed nature of the pointing directions \cite{Sra:Axial:2013}. With this, directional statistics can be applied once again.

\subsection{Regression Revisited}
\label{subsec: Regression Revisited}

Let us apply the averages developed in this section to the equations in section~\ref{sec: application to regression}, \emph{Application to Regression}. The \SVD and orientation steps remain the same. The linear regression of~(\ref{eq: regression with oriented basis}) is modified to use the local average of the eigenvalues from~(\ref{eq: eigenvalue average}),
\begin{equation}
  y = U S\: \bar{\EigValMatrix}^{1/2}\: \hat{\beta} + \epsilon.
  \label{eq: regression with oriented basis and avg e-values}
\end{equation}
The prediction step is also modified from~(\ref{eq: E y = P R beta}) to read
\begin{equation}
  \Expectation y = P_{\text{out}}\: \bar{\orientedeigvec}\; \hat{\beta}.
  \label{eq: E y = R V-bar beta}
\end{equation}
Here the mean direction of the eigenvectors in~(\ref{eq: mean direction V-bar def}) replaces the single-observation rotation matrix~$R$ in the original. 

Use of eigenvalue and eigenvector averages will reduce sample-based fluctuation and therefore may improve the predictions. However, from a time-series perspective, averages impart delay because an average is constructed based on a lookback interval. For stationary underlying processes the delay may not a problem, yet for nonstationary processes the lag between the average and the current state can lead to lower quality predictions. 

In any event, at the very least this section has presented a methodology to decompose variation across a sequence of eigensystem observations into meaningful quantities.

%
%

\section{Reference Implementation}
\label{sec: reference implementation}

The Python package \code{thucyd}, written by this author, is freely available from \textsc{PyPi}\footnote{%
Available at \href{https://pypi.org/project/thucyd/}{\code{pypi.org/project/thucyd}}
}
and \textsc{Conda-Forge}\footnote{%
Available at \href{https://github.com/conda-forge/thucyd-feedstock}{\code{github.com/conda-forge/thucyd-feedstock}}
}
and can be used directly once installed. The source code is available on at \href{https://gitlab.com/thucyd-dev/thucyd}{\code{gitlab.com/thucyd-dev/thucyd}}. 

There are two functions exposed on the interface of \code{thucyd.eigen}:
\begin{itemize}
  \item \code{orient\_eigenvectors} implements the algorithm in section~\ref{sec: algorithm details}, and
  \item \code{generate\_oriented\_eigenvectors} implements \\  $R_k = \mathcal{G}(\boldsymbol{\theta}, k)$, eq. (\ref{eq: R = G(theta) def}).
\end{itemize}
The pseudocode in listing~\ref{eigen::eq: algorithm pseudocode} outlines the reference implementation for \code{orient\_} \code{eigenvectors} that is available at the above-cited source-code repositories. Matrix and vector indexing follows the Python Numpy notation.

\SetAlCapSty{}
\begin{algorithm}[t]
  {
    \DontPrintSemicolon
    \SetKw{KwStep}{step}
    \SetKwData{V}{V}
    \SetKwData{Vcol}{Vcol}
    \SetKwData{Vor}{Vor}
    \SetKwData{Vsort}{Vsort}
    \SetKwData{Vwork}{Vwork}
    \SetKwData{E}{E}
    \SetKwData{Eor}{Eor}
    \SetKwData{Esort}{Esort}
    \SetKwData{R}{R}
    \SetKwData{cursor}{cursor}
    \SetKwData{cursors}{cursors}
    \SetKwData{SignFlipVec}{SignFlipVec}
    \SetKwData{AnglesCol}{AnglesCol}
    \SetKwData{AnglesMtx}{AnglesMtx}
    \SetKwData{AnglesWork}{AnglesWork}
    \SetKwData{SortIndices}{SortIndices}
    \SetKwFunction{OrientEigenvectors}{OrientEigenvectors}
    \SetKwFunction{SortEigenvectors}{SortEigenvectors}
    \SetKwFunction{ReduceDimensionByOne}{ReduceDimensionByOne}
    \SetKwFunction{SolveRotationAnglesInSubDim}{SolveRotationAnglesInSubDim}
    \SetKwFunction{ConstructSubspaceRotationMtx}{ConstructSubspaceRotationMtx}
    \SetKwFunction{MakeGivensRotation}{MakeGivensRotation}
    \SetKwProg{Def}{def}{:}{}
    %
    \Def{\OrientEigenvectors{\V, \E}}{
      \Vsort, \Esort, \SortIndices $\leftarrow$ \SortEigenvectors{\V, \E}\;
      \Vwork $\leftarrow$ copy(\Vsort), $N$ $\leftarrow$ $\text{dim}(\V)$ \;
      \For{$i=0$ \KwTo $N-1$}{
        $\SignFlipVec[i]$ $\leftarrow$ $(\Vwork[i,i] => 0)$ ? $1$ : $-1$ \;
        $\Vwork[:, i]$ $\leftarrow$ $\Vwork[:, i]\cdot\SignFlipVec[i]$\;
        \Vwork, $\AnglesMtx[i, :]$ $\leftarrow$ \ReduceDimensionByOne{$N$, $i$, \Vwork} \;
      }
      \Vor $\leftarrow$ $\Vsort\cdot\diagoper(\SignFlipVec)$, \Eor $\leftarrow$ \Esort \;
      return \Vor, \Eor, \AnglesMtx, \SignFlipVec, \SortIndices   \;
    }
    \;
    %
    \Def{\SortEigenvectors{\V, \E}}{
      \SortIndices $\leftarrow$ argsort$(\text{abs}(\diagoper(\E)))$ \;
      return $\Vsort[:, \SortIndices]$, $\diagoper(\Esort[\SortIndices])$, \SortIndices \;
     }
    \;
    %
    \Def{\ReduceDimensionByOne{$N$, $i$, \Vwork}}{
      \Vcol $\leftarrow$ $\Vwork[:, i]$ \;
      \AnglesCol $\leftarrow$ \SolveRotationAnglesInSubDim{$N$, $i$, \Vcol} \;
      \R $\leftarrow$ \ConstructSubspaceRotationMtx{$N$, $i$, \AnglesCol} \;
      \Vwork $\leftarrow$ $\R\transpose\Vwork$ \;
      return \Vwork, $\AnglesCol\transpose$ \;
    }
    \;
    %
    \Def{\SolveRotationAnglesInSubDim{$N$, $i$, \Vcol}}{
      \AnglesWork $\leftarrow$ zeros$(0: N)$ \;
      r $\leftarrow$ 1 \;
      \For{ $j = N$  \KwTo $i + 1$ \KwStep $-1$}{
        $y$ $\leftarrow$ $\Vcol[j]$  \;
        $r$ $\leftarrow$ $r * \text{cos}(\AnglesWork[j+1])$ \;
        $\AnglesWork[j]$ $\leftarrow$ $(r\; != 0)$ ? arcsin$(y / r)$ : 0 \; 
      }
      return $\AnglesWork[:N]$ \;
    }   
    \;
    %
    \Def{\ConstructSubspaceRotationMtx{$N$, $i$, \AnglesCol}}{
      \R $\leftarrow$ $I(N)$ \;
      \For{ $j = N$ \KwTo $i+1$ \KwStep $-1$ }{
         \R $\leftarrow$ $\MakeGivensRotation{$N, i, j, \AnglesCol[j]$} \cdot \R$  \;
      }
      return \R \;
    }
    \;
    %
    \Def{\MakeGivensRotation{$N, i, j, \theta$}}{
      \R $\leftarrow$ $I(N)$, $c$ $\leftarrow$ $\cos(\theta)$, $s$ $\leftarrow$ $\sin(\theta)$ \;
      \R$[i,i]$ $\leftarrow$ $c$, \R$[j,j]$ $\leftarrow$ $c$, 
      \R$[i,j]$ $\leftarrow$ $s$, \R$[j,i]$ $\leftarrow$ $s$ \;
      return \R \;
    }
    \;
  }
\caption{Pseudocode implementation of \code{orient\_eigenvectors}.}
\label{eigen::eq: algorithm pseudocode}
\end{algorithm}

\subsection{Eigenvector Orientation}
\label{subsec: eigenvector orientation}

A simple orientation example in~$\realsspace^3$ reads as follows:
\begin{lstlisting}[language=Python,caption=]
import numpy as np
import thucyd

# setup an example eigenvector mtx with qr
full_dim = 3; A = np.eye(full_dim)
A[:, 0] = 1. / np.sqrt(full_dim)
V, _ = np.linalg.qr(A)
# setup trivial eigenvalues, descending order
E = np.diag(np.arange(full_dim)[::-1])
# reflect the second eigenvector to build Vin
S = np.diag(np.array([1., -1., 1.]))
Vin = V.dot(S)
# orient the original matrix V
Vor, Eor, signs, theta_matrix, sort_indices = \
    thucyd.eigen.orient_eigenvectors(Vin, E)
\end{lstlisting}
Execution of this snippet\footnote{%
The snippet was run with Numpy version 1.16.4 on MacOS 10.14.5. OpenBlas and Intel MKL produced the same results, as expected.
}
converts the original basis 
\begin{equation*}
	\begin{aligned}
	  \code{Vin} &= \left(
	    \begin{array}{ccc}
	      -0.577 & -0.408 & -0.707 \\
	      -0.577 & +0.816 &  0     \\
	      -0.577 & -0.408 & +0.707 
	    \end{array}
	  \right)
	  \quad\text{to}\\
	  \code{Vor} &= \left(
	    \begin{array}{ccc}
	        0.577 & -0.408 & -0.707 \\
	        0.577 & +0.816 &  0     \\
	        0.577 & -0.408 & +0.707 
	    \end{array}
	  \right),
	\end{aligned}
\end{equation*}
with \ebox{\code{signs} = \left(-1, +1, +1\right)}. The cross product exists in~$\realsspace^3$ so the column vectors in \code{Vin} can be inspected to determine their handedness. In this case, \code{Vin} is a left-handed basis, so a pure rotation cannot align the basis of \code{Vin} to~$\Identity$. However, once the first eigenvector is reflected, a pure rotation is all that is required for alignment. The matrix of rotation angles \code{theta\_mtx}, which follows~(\ref{eq: theta mtx storage}) in form and is expressed here in degrees, is
\begin{equation*}
  \boldsymbol{\theta} = \left(
    \begin{array}{ccc}
         0  &  \;45\;  &   35.264  \\
         0  &  0   &  -30      \\
         0  &  0   &   0
    \end{array}
  \right).
\end{equation*}
The generator function consumes this array of angles to reconstruct the oriented eigenbasis, as in
\begin{lstlisting}[language=Python,caption=]
# reconstruct Vor from rotation angle matrix
Vor_recon = thucyd.eigen.generate_oriented_eigenvectors(theta_mtx)
\end{lstlisting}
Running this snippet will show that \ebox{\code{Vor\_recon} = \code{Vor}}.

\subsection{Eigenvector Reconstruction from Rotations}
\label{subsec: eigenvector reconstruction from rotations}

The reconstruction of an oriented eigenvector from rotations is better illustrated in~$\realsspace^4$. Here is a snippet that calls the \code{generate\_oriented\_eigenvectors} api function.
\begin{lstlisting}[language=Python,caption=]
# setup example eigenvector and eigenvalue mtx
full_dim = 4; A = np.eye(full_dim)
A[:, 0] = 1. / np.sqrt(full_dim)
Vin, _ = np.linalg.qr(A)
E = np.diag(np.arange(full_dim)[::-1])
# cast Vin into an oriented basis
Vor, _, _, theta_matrix, _ = impl.orient_eigenvectors(Vin, E)
# build W_k
W_k = np.eye(full_dim)
# iterate across columns in (Vor, W_k)
for cursor in np.arange(full_dim):
    # call focus fcxn (the 2nd argument is in base(1))
    W_k = W_k.dot(
        thucyd.eigen.generate_oriented_eigenvectors(
            theta_matrix, cursor + 1))
# reconstructed eigenbasis
Vor_recon = W_k
\end{lstlisting}
The output \code{theta\_matrix} from the function \code{orient\_} \code{eigenvectors} is the input to the focus function here. While we could simply call
\begin{lstlisting}[language=Python,caption=]
Vor_recon = thucyd.eigen.generate_oriented_eigenvectors(theta_matrix)
\end{lstlisting}
it is revealing to build up \code{Vor\_recon} one subspace at a time. 

The oriented eigenvector and angles matrices are
\begin{equation*}
  \begin{aligned}  
	  \code{Vor} &= \left(
	    \begin{array}{cccc}
	        0.5   & -0.289 & -0.408 & -0.707    \\
	        0.5   & +0.866 &  0.000 &  0.000    \\
	        0.5   & -0.289 & +0.816 &  0.000    \\
	        0.5   & -0.289 & -0.408 &  0.707
	    \end{array}
	  \right)
	  \quad\text{and}\\
	  \boldsymbol{\theta} &\simeq \left(
	    \begin{array}{cccc}
	        0    &  45   &   35 &  30       \\
	        0    &   0   &  -30     & -19   \\
	        0    &   0   &    0     & -30   \\
	        0    &   0   &    0     &  0
	    \end{array}
	  \right).
  \end{aligned}
\end{equation*}
We can now build up \code{Vor} through \code{W\_k} from these embedded angles:
\begingroup
\allowdisplaybreaks
  \begin{align*} 
    R_1 & = 
        \left(
        \begin{array}{cccc}
            0.5   & -0.707 & -0.408 & -0.289  \\
            0.5   &  0.707 & -0.408 & -0.289  \\
            0.5   &  0.000 &  0.816 & -0.289  \\
            0.5   &  0.000 &  0.000 &  0.866  
        \end{array}
        \right) \;, \\
    R_1 R_2 & = 
        \left(
        \begin{array}{cccc}
            0.5   & -0.289 & -0.707 & -0.408  \\
            0.5   & +0.866 &  0.000 &  0.000  \\
            0.5   & -0.289 &  0.707 & -0.408  \\
            0.5   & -0.289 &  0.000 &  0.816  
        \end{array}
        \right) \;, \\
    R_1 R_2 R_3 & = 
        \left(
        \begin{array}{cccc}
            0.5   & -0.289 & -0.408 & -0.707  \\
            0.5   & +0.866 &  0.000 &  0.000  \\
            0.5   & -0.289 & +0.816 &  0.000  \\
            0.5   & -0.289 & -0.408 &  0.707  
        \end{array}
        \right) \;, \quad\text{and} \\
    R_1 R_2 R_3 R_4 & = 
        \left(
        \begin{array}{cccc}
            0.5   & -0.289 & -0.408 & -0.707  \\
            0.5   & +0.866 &  0.000 &  0.000  \\
            0.5   & -0.289 & +0.816 &  0.000  \\
            0.5   & -0.289 & -0.408 &  0.707  
        \end{array}
        \right) \;.
  \end{align*}
\endgroup
It is apparent that as the subspace rotations are concatenated, the rightward columns become aligned to the reference eigenbasis \code{Vor}. That the reconstructed matrix matches the reference matrix before the last rotation is simply because the last rotation, $R_4$, is the identity matrix, and is only included in the formalism for symmetry. 
 
Looking ahead, an optimized implementation can avoid trigonometry by carrying the sine functions \ebox{\sin\theta} rather than the angle values~$\theta$. On the range \ebox{\theta\in[-\pi/2,} \ebox{ \pi/2]} the function is invertible, thus either representation will do. Elimination of trigonometry leaves only arithmetic and matrix multiplication, thus making the \code{orient} algorithm a candidate for \LAPACK implementation.

%
%

\section{Conclusions}
\label{sec: conclusions}

Although eigenanalysis is an old and well-studied topic, linking eigenanalysis to an evolving dataset leads to unexpected results and new opportunities for analysis. The optimization of eigenanalysis codes for one-time solutions has admitted inconsistent eigenvector orientation because such inconsistency is irrelevant to the one-time solution. Yet for an evolving system it is precisely the inconsistency that disrupts interpretability of an eigen-based model. This article reports a method to correct for the inconsistency, assuming a framework where orientation is a postprocessing step to existing eigenanalysis implementations. Once corrected, directional statistics brings a new avenue of inquiry to the behavior of the underlying factors in the dataset.

It is hoped that the methods in this paper will be applied in fields as disparate as human health, climate change, navigation systems, economics, the internet, and other scientific fields.


\begin{acknowledgements}

Thank you to Professor A. Brace (National Australia Bank, University of Technology Sydney), %
Professor. P. N. Kolm (New York University), Dr. A. Meucci (Advanced Risk and Portfolio Management), %
Professor H. R. Miller (Massachusetts Institute of Technology), and %
Dr. J. P. B. M{\"u}ller (Rhizome Works) %
for their helpful conversations about, and encouragement of, this publication. 

\end{acknowledgements}

%
%




\end{document}